
\magnification=\magstep1
\vsize=22truecm
\input amstex

\input epsf.tex

\TagsOnRight

\def\cite#1{{\rm [#1]}}
\def\ol{\overline }




\document

{\centerline{\bf{GENERALIZED FORMS OF AN OVERCONSTRAINED}}}
{\centerline{\bf{SLIDING MECHANISM CONSISTING}}}
{\centerline{\bf{OF TWO CONGRUENT TETRAHEDRA}}}

\medskip

\medskip

{\centerline{ENDRE MAKAI, JR.$^*$, TIBOR TARNAI$^{**}$}}
\vskip1.0cm

Endre Makai, Jr., 
Alfr\'ed R\'enyi Mathematical Institute, 
L. E\"otv\"os Research Network (ELKH), 
H-1364 Budapest, P.O. Box 127, 
HUNGARY,
makai.endre\@renyi.hu,
\newline
{\rm{http://www.renyi.hu/\~{}makai}},
ORCID ID: https://orcid.org/0000-0002-1423-8613

Tibor Tarnai, 
Budapest University of Technology and Economics, 
Department of Structural Mechanics, 
H-1521 Budapest, M\H{u}egyetem rkp.\ 3, 
HUNGARY,
\newline
tarnai.tibor\@emk.bme.hu, {\rm{http://www.me.bme.hu/tarnai-tibor,}}
\newline
ORCID ID: https://orcid.org/0000-0001-9260-7800

\medskip

{\it 2020 Mathematics Subject Classification:} 51M99
\vskip.0cm
{\it Keywords and phrases:} {\rm{Sliding mechanisms,
tetrahedra, bar structures}}

\medskip

$^*$Research (partially) supported by Hungarian National Foundation for
Scientific Research, grant nos.\ K68398, T046846.\newline
\indent $^{**}$Research (partially) supported by Hungarian National Foundation
for
Scientific Research, grant no.\ T046846 and NKFI, grant
no.\ K138615.



\medskip

ABSTRACT. We investigate the
motions of a bar structure consisting of two congruent
tetrahedra, whose edges in their basic position form the face diagonals of a
rectangular parallelepiped.
The constraint of the motion is that the originally intersecting edges should
remain coplanar.
We determine all finite motions of our bar structure.
This generalizes our earlier work, where we did the same for the case when the
rectangular parallelepiped was a cube. At the end of the
paper we point out three further possibilities to generalize 
the question about the cube, and give for them examples of
finite motions.



\medskip

\medskip

{\bf{1.}} INTRODUCTION

\medskip

Drawing all the diagonals of all faces of a cube, we obtain the edges of two
congruent regular tetrahedra.
This position of these tetrahedra is called their
{\it basic position}.

\smallskip
{\narrower{\narrower\noindent
We keep one of the tetrahedra fixed. We move the other one
under the following condition.
Each pair of edges of the two tetrahedra, which were
originally diagonals of
some face of the cube, should remain coplanar.
\par}}

\smallskip
\noindent
The structure consisting
of the above described two tetrahedra has been invented in 1982
by L. Tompos, Jr. He was then an undergraduate of the Hungarian Academy of Craft and
Design.
He has built a physical model of the bar (i.e., edge)
structure of these tetrahedra, as follows. 
The bars of one of the tetrahedra touched those of the other tetrahedron
from inside (Fig.~1). He has observed

\noindent
FIGURE 1 ABOUT HERE
%
\smallskip
\vbox{
\centerline{\epsfxsize=61mm 
\epsfbox{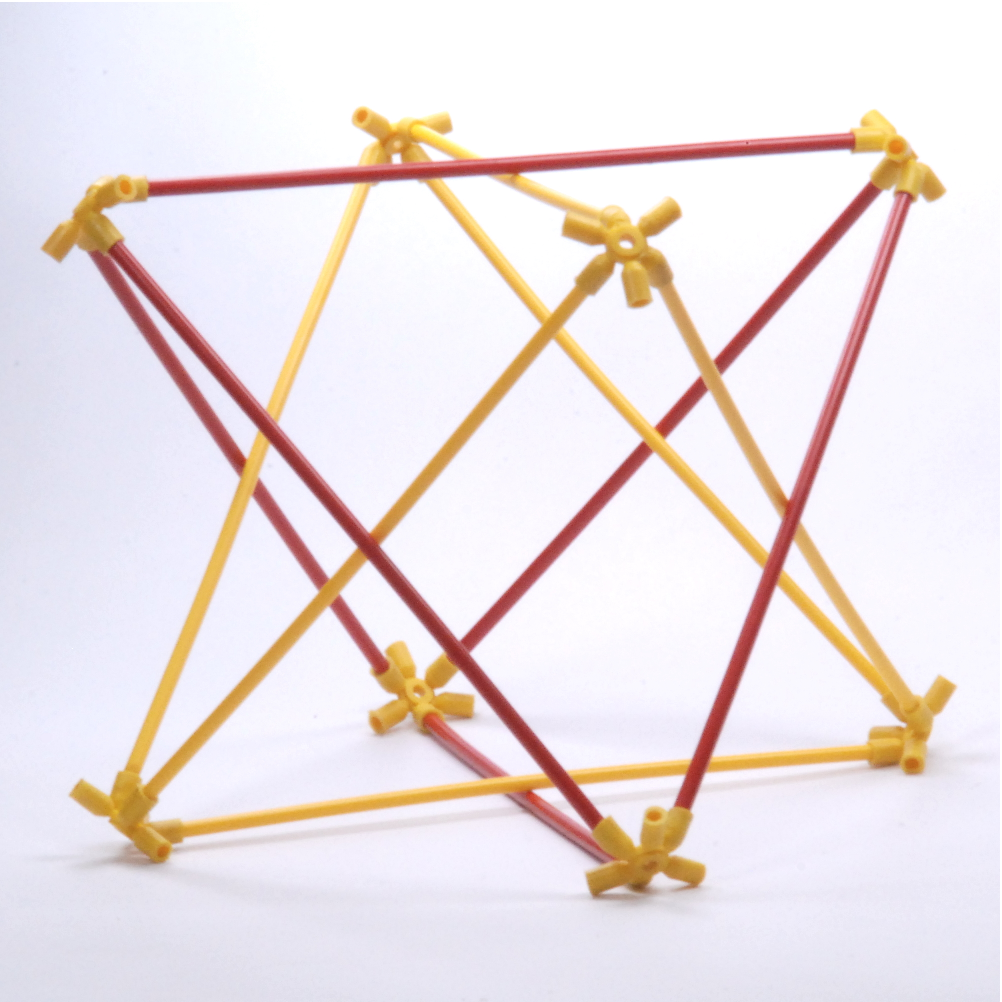}}
\smallskip
\noindent
{
{\centerline{Figure 1.
The bar-and-joint structure of Tompos's pair of}}
{\centerline{
tetrahedra. The physical model in the basic position.}}
\newline
{\centerline{(Photograph provided by Andr\'as Lengyel.)}}
\par}
}
\smallskip
%
%
\noindent
that this structure admits continuous motions.
We note that \cite{2},
p.\ 7 contains a figure of these tetrahedra, but
their mobility is not investigated there. We have found the
same figure as a decoration of the dining room of 
Hotel Arcas de Agua
in Spain, in village Arcas, Cuenca (Fig.~2).

\noindent
FIGURE 2 ABOUT HERE
%
\smallskip
\vbox{
\centerline{\epsfxsize=61mm 
\epsfbox{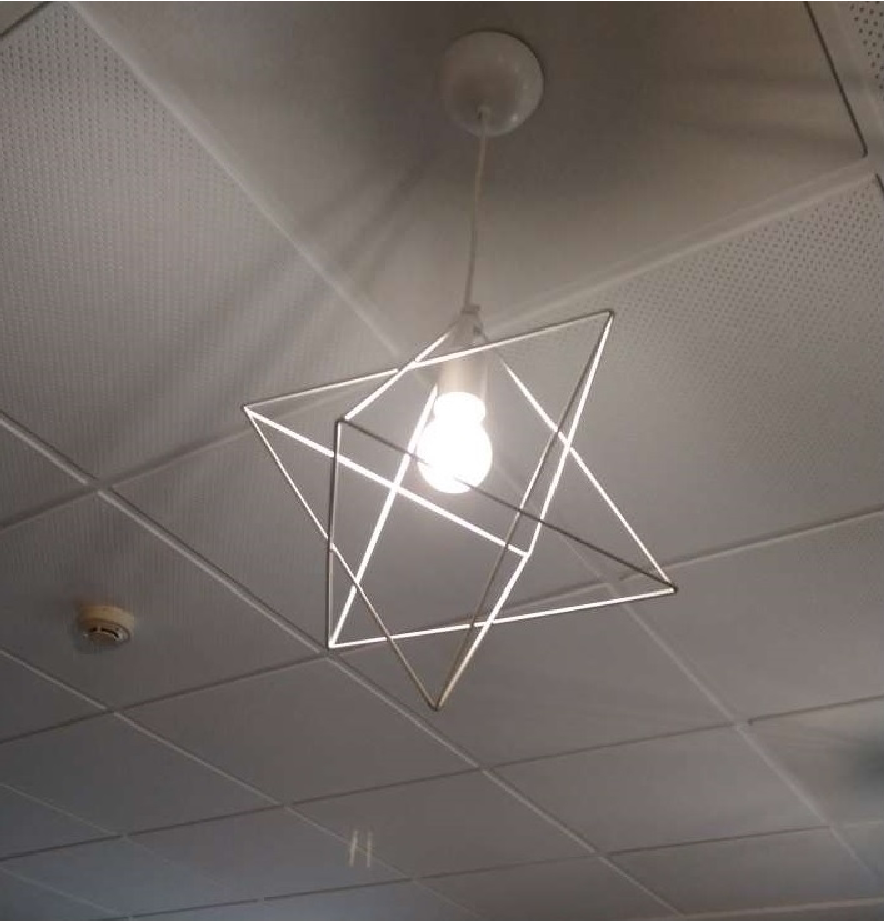}}
\smallskip
\noindent
{
{\centerline{Figure 2.
The lamp decoration in the hotel in Arcas.}}
{\centerline{(Photograph provided by Ampar L\'opez.)}}
\par}
}
\smallskip
%

By a {\it motion\/} (sometimes we will say a
{\it finite motion}) we will not mean
a continuous motion from the basic position,
always satisfying the constraints,
but {\it any position of our
structure that satisfies the constraints\/}. (Possibly this position is not the
result of a continuous motion,
always satisfying the constraints.)
We mean by this the following. One of the tetrahedra is fixed.
The other one
is obtained
from the basic position of itself by the application of an isometry
(i.e.,
congruence) of the space {\it of determinant~$+1$}. Additionally,
the coplanarity
conditions are satisfied.
An isometry of determinant~$+1$ will be written in the form
$\bold{\boldsymbol\Phi(x) = A x + b}$. Here $\bold A$ is a $3 \times 3$
orthogonal matrix of determinant~$+1$, and $\bold b$ is a vector in~$\bold R^3$.
Geometrically, $\bold A$ is a rotation about some straight line containing the
origin.

\cite{10} and \cite{12} determined all motions
of this pair of tetrahedra.
We give a brief description of them, with the names used for them
in \cite{12}.
Note that the bars of the model have non-zero width. Hence
not all motions can be
realized by the physical model. We can realize
only such motions for which each respective pair of
edges (which have to be coplanar) actually has a common point. 
A discussion of the question
which motions are physically admissible, i.e., satisfy this
more restrictive condition, is contained in \cite{1}, Ch.~4.
Further in this paper we will make no distinction between
physically
admissible and inadmissible motions.

We suppose that the vertices of the tetrahedra in the basic position are the
points $(\pm 1, \pm 1, \pm 1)$.
More exactly, the vertices of the fixed tetrahedron are $P^0_1(1, -1, -1)$,
\newline
$P^0_2(-1, 1, -1)$, $P^0_3(-1, -1, 1)$, $P^0_4(1, 1, 1)$.
The vertices of
the moving tetrahedron are $Q^0_1, \dots , Q^0_4$.
Here, in the basic
position, $Q^0_i$ is the mirror image of $P^0_i$ w.r.t.\ the origin (i.e., the
centre of the cube with vertices $(\pm 1, \pm 1, \pm 1)$), see Fig.~3.

\noindent
FIGURE 3 ABOUT HERE
%
\smallskip
\vbox{
\centerline{\epsfxsize=86mm 
\epsfbox{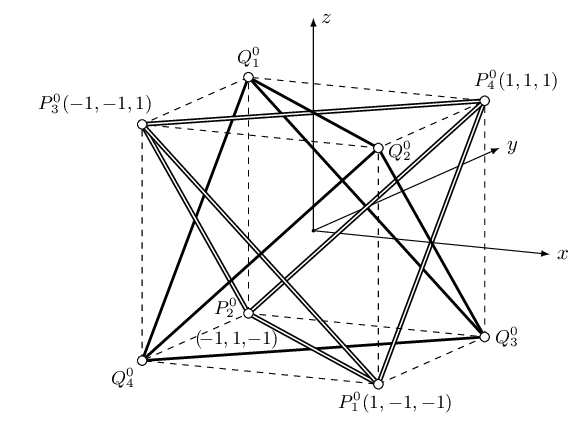}}
\smallskip
\noindent
{
{\centerline{Figure 3.
The notation for the vertices of the tetrahedra.}}
\par}
}
\smallskip
%

There exist motions $\bold{\boldsymbol\Phi(x) = Ax + b}$ of
Tompos's
tetrahedra, where $\bold A$ is a rotation about an axis $\bold{0 e}_i$, $\bold
0(\bold e_i \pm \bold e_j)$ or $\bold 0(\bold e_i \pm \bold e_j \pm \bold
e_k)$, resp. Here $\bold e_1, \bold e_2, \bold e_3$ are the basic unit
vectors in the space, and $i, j, k$ are different.
These motions $\bold{\boldsymbol\Phi(x)}$ are called
the {\it{motions of the first,
second and third kind}}, resp. In the first case, for angle of rotation
$\pi $, ${\bold{b}}$ is not unique, but
we count only the case $\bold{b = 0}$ to the motion of the first
kind.
The angle of rotation of $\bold A$ in the first case is arbitrary, in the
second case it is arbitrary, except $\pi $, in the third case it is
arbitrary, except $\pm \pi /2$.
(The angle of rotation is positive if, looking from the axis
vector, e.g., ${\bold{0e}}_i$, backwards,
it is positive.)

It turned out that for each above rotation $\bold A$, except in the first case
the rotation through the angle
$\pi $, we had the following. There existed a unique translation $\bold b$
such that $\bold{\boldsymbol\Phi(x) = Ax + b}$ was a motion of
Tompos's
tetrahedra.
In the first case, for an angle of rotation different from $\pi $, we
had $\bold{b = 0}$.

There also exist motions $\bold{\boldsymbol\Phi(x) = A x + b}$ of Tompos's
tetrahedra, with $\bold A$ a rotation about an axis $\bold{0}(C_1 \bold e_i +
C_2 \bold e_j)$, where $i \neq j$, and $C_1, C_2$ are real, not both~$0$.
If $C_1 C_2 \neq 0$, and so this is not a motion of the first kind, then the
angle of the rotation $\bold A$ is arbitrary, except $\pi $. Moreover, for
each such rotation $\bold A$ the translation $\bold b$ is uniquely determined.
These motions, for $C_1 C_2 \neq 0$, together with the motions of the first
kind (where $C_1 C_2 = 0$), both with angle of rotation
different from $\pi$, are called the
{\it{motions of the intermediate kind}}.
These motions also contain the motions of the second kind as a special case.

Now let $\bold A$ be a rotation about an axis $\bold 0 \bold e_i$, through the
angle
$\pi $. Then let $\bold b$ be any vector of the form $C\bold e_j$, or
$C_1 \bold e_j + C_2 \bold e_k$, resp., where $C, C_1, C_2$ are real
and $i, j, k$ are different.
Then $\bold{\boldsymbol\Phi(x) = Ax + b}$ is a motion of Tompos's
tetrahedra, which is called the
{\it{motion of the fourth, or fifth kind}},
resp.
The motions of the fifth kind contain the motions of the fourth kind as a
special case.
The motions of the first, second, third, fourth and fifth kinds are drawn in
Figs.~4a, 4b, 4c, 4d, 4e.

\noindent
FIGURES 4a, 4b, 4c, 4d, 4e ABOUT HERE 
%
\smallskip
\vbox{
\centerline{\epsfxsize=111mm 
\epsfbox{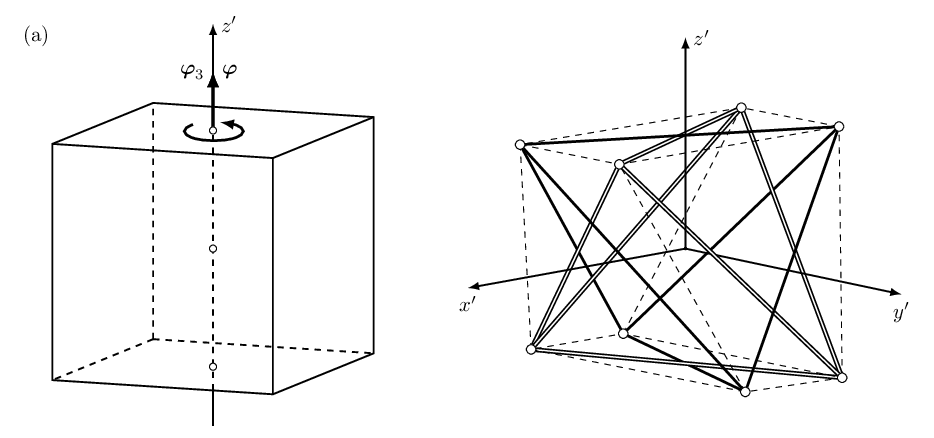}}
\smallskip
\noindent
{
{\centerline{Figure 4a.
The motion of the first kind.}}
\par}
}
\smallskip
%

%
\smallskip
\vbox{
\centerline{\epsfxsize=111mm 
\epsfbox{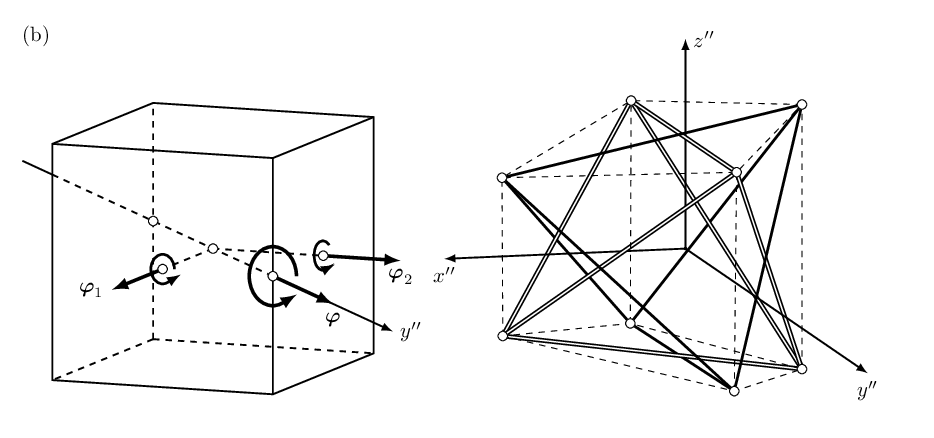}}
\smallskip
\noindent
{
{\centerline{Figure 4b.
The motion of the second kind.}}
\par}
}
\smallskip
%

%
\smallskip
\vbox{
\centerline{\epsfxsize=111mm 
\epsfbox{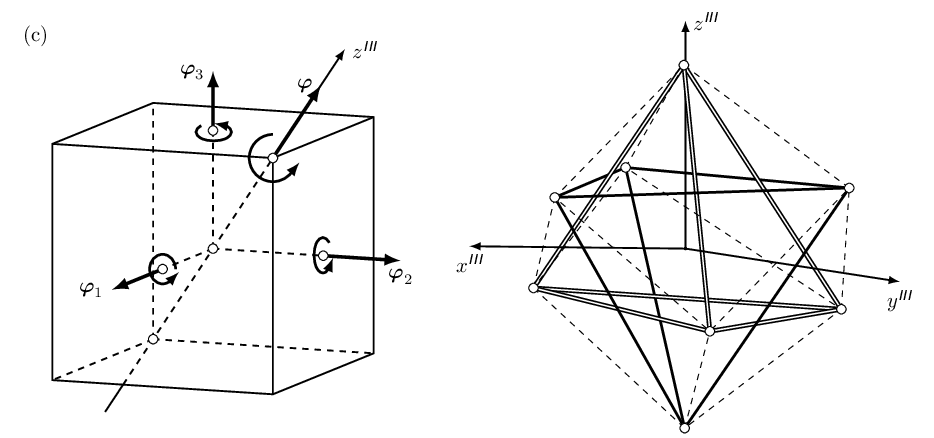}}
\smallskip
\noindent
{
{\centerline{Figure 4c.
The motion of the third kind.}}
\par}
}
\smallskip
%

%
\smallskip
\vbox{
\centerline{\epsfxsize=111mm 
\epsfbox{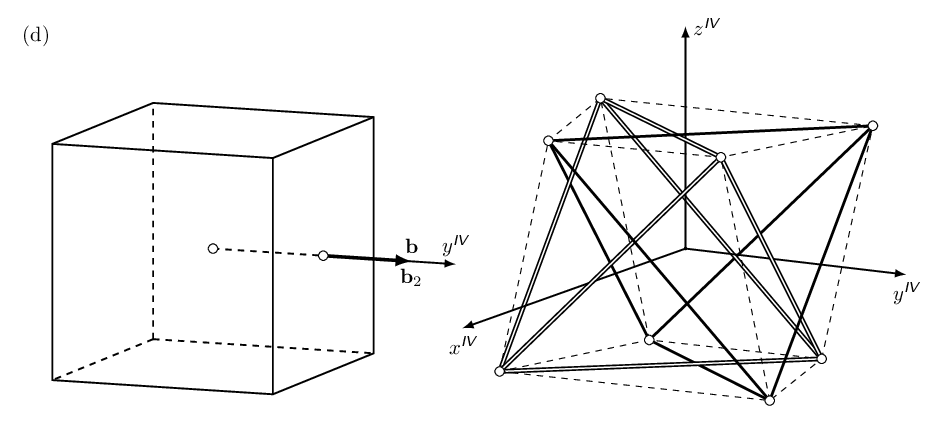}}
\smallskip
\noindent
{
{\centerline{Figure 4d.
The motion of the fourth kind.}}
\par}
}
\smallskip
%

%
\smallskip
\vbox{
\centerline{\epsfxsize=111mm 
\epsfbox{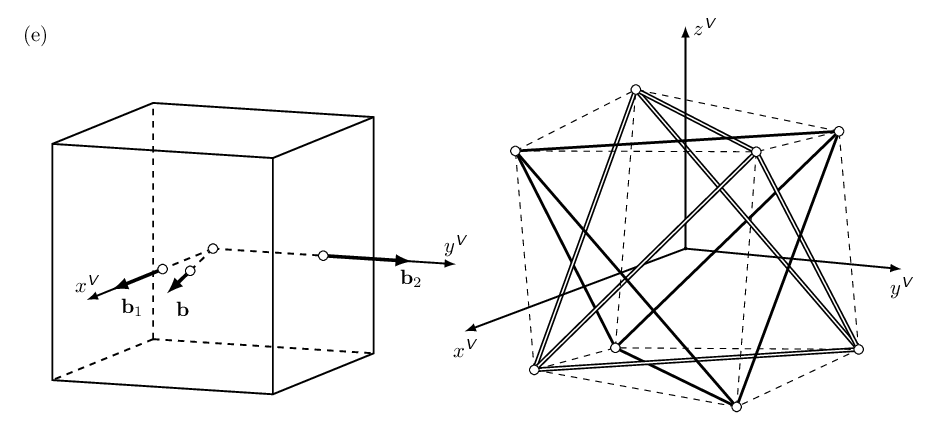}}
\smallskip
\noindent
{
{\centerline{Figure 4e.
The motion of the fifth kind.}}
\par}
}
\smallskip
%

It turned out that the motions of each kind, for given
${\bold{e}}_i$,
${\bold{e}}_i \pm {\bold{e}}_j$, ${\bold{e}}_i \pm
{\bold{e}}_j \pm {\bold{e}}_k$, $\{ {\bold{e}}_i,
{\bold{e}}_j \} $,
${\bold{e}}_j$ and $\{ {\bold{e}}_j, {\bold{e}}_k \} $
(in the above order), resp.,
constituted a smooth manifold in the six-dimensional
manifold of
all motions (i.e., isometries of determinant $+1$)
of the space. They had dimensions $1, 1, 1, 2, 1, 2$,
resp.\ (in the above order).
The motions of the third kind formed a manifold of two connected
components,
cf.\ \cite{13}, p.~141.

These manifolds show certain bifurcation phenomena, which have been analyzed
in \cite{12}.
Moreover, as shown by \cite{10} and \cite{12}, the
above enumerated motions are the only motions of Tompos's tetrahedra.
We note that \cite{10} also described the trajectories of the
vertices during the physically admissible motions.
Moreover, in
\cite{4, 5, 6} (of which
only \cite{6} has been available to the authors)
and in \cite{1}
the motions of Tompos's tetrahedra have been further
investigated. They have pointed out
some possible mechanical engineering applications.

In our paper we generalize the above investigations.
We start not with a cube, but with a general
rectangular parallelepiped.
All diagonals of all of its faces constitute the edges of two congruent
tetrahedra.
This position of these two tetrahedra is called their
{\it basic position}.

\smallskip
{\narrower{\narrower\noindent
We keep one of the tetrahedra fixed. We move the other one
(i.e., apply to it
an isometry of the space, of determinant~$+1$) under the
following condition.
Each pair of edges of the two tetrahedra, which
were originally diagonals
of some face of the rectangular parallelepiped, should
remain coplanar.
\par}}

\vskip-12.5pt
\hfill (A)

\smallskip
\noindent
In the physical model, the bars (edges) of one of the
tetrahedra touch those of the other tetrahedron from inside,
as in Fig.\ 1.
First we give the description of all such
motions, and then we give a mathematical
proof that this list of motions is complete.
The results are rather analogous to the case of the cube. Only in a special
case there is a motion of the sixth kind, which constitutes a
$1$-manifold
(one-dimensional
manifold). These results, more exactly, each of Theorems 1,
2 and 3 of this paper, have been announced in
\cite{8}. We also analyze the bifurcation properties of
the solution manifolds.

Second, we give three further generalizations of our pair of
tetrahedra. For these we will not be able to determine all finite
motions, but we will be able to give certain finite motions.
These generalizations
are the following.
(1)
A pair of tetrahedra derived from a general parallelepiped.
(2)
A pair of regular $n$-gonal pyramidal frames (the
bases are allowed to change their shape, and also
to become non-planar).
(3)
A pair of regular tetrahedra with congruent circular arc
edges.
In case (2)
we will present some numerical evidence that
all continuous finite motions from the basic position,
always satisfying the constraints, 
might form the one-parameter family of
finite motions which we have found. Some of these results,
namely concerning (1) and (2),
also have been announced in \cite{8}.

Related results on pairs of polyhedra or polyhedral frames
moving with sliding
constraints cf.\ in \cite{7} and \cite{9}.


\medskip

\medskip

{\bf{2.}} THE MOTIONS OF THE TWO TETRAHEDRA DERIVED FROM A
RECTANGULAR PARALLELEPIPED

\medskip

{\bf 2.1.}
We will use analogous notations as in the case of the cube.
Let the vertices of the rectangular parallelepiped be
$(\pm d_1, \pm d_2, \pm
d_3)$, where $d_1, d_2, d_3 > 0$ are constant.
The fixed vertices are $P_1(d_1, -d_2, -d_3)$, $P_2(-d_1, d_2, -d_3)$,
$P_3(-d_1, -d_2, d_3)$, $P_4(d_1, d_2, d_3)$. 
The moving vertices are
$Q_1$, $Q_2$, $Q_3$, $Q_4$, where in the basic position $Q_i$ is the mirror
image of $P_i$ w.r.t.\ the origin.
Thus $P_1 P_2 P_3 P_4$ is the fixed tetrahedron, and
$Q_1 Q_2 Q_3 Q_4$ is the
moving tetrahedron. We move the moving tetrahedron under
condition (A), which is the same as
in the case of the cube.
First we describe the motions $\bold{\boldsymbol\Phi(x) = Ax + b}$ of this
moving tetrahedron.

First of all, the motions of the
fourth and fifth kinds are defined word for word
as in the case of the cube. They evidently exist (i.e., the originally
intersecting edges remain coplanar).

The motion of the first kind, also  defined in the same way
as for the cube,
exists for any value of the rotation angle.
Moreover, we have for this motion $\bold{b = 0}$.
These are shown as in
\cite{12}, p.~428.
Namely, let the axis of rotation be, e.g.,
${\bold{0}}{\bold{e}}_3$. Then the edges of the moving
tetrahedron, originally lying on some horizontal face of the
rectangular parallelepiped, contain $\pm d_3{\bold{e}}_3$.
Hence they remain coplanar. Now consider
the straight lines spanned by the edges of the moving
tetrahedron, originally lying on some
two opposite vertical faces of the rectangular
parallelepiped. Their rotated copies, through all angles,
constitute a ruling of a (doubly) ruled surface, a one-sheet
hyperboloid of revolution.
The other diagonals of these two faces belong
to the other ruling of this hyperboloid. However, any two
lines from different rulings of this hyperboloid are
coplanar. Hence this motion exists for all $\varphi $.

The motion of the intermediate kind is also defined in
the same way as for the cube.
(Now we do not need to treat separately any analogue of the
motion of the second kind.)
The motion of the intermediate kind exists for
any value of the rotation angle, different from $\pi $. (But
sometimes it may exist also for rotation angle $\pi $, cf.\
below.)
We apply the analogues of
the arguments in 
\cite{12}, p.\ 429.
Let the axis of rotation be ${\bold{0u}}$, where
${\bold{u}} = C_1{\bold{e}}_1 + C_2{\bold{e}}_2$ has length
$1$.

As in \cite{12},
for convenience, rather than taking one tetrahedron as
fixed, the other one as
moving, we make the following.
We rotate both tetrahedra about the rotation axis,
through angles $\pm \varphi / 2$, in such a way that their
symmetry w.r.t.\ the $xy$-plane is preserved
(thus they are rotated in opposite senses). Then the
edges of the two tetrahedra, originally lying on one
vertical face of the rectangular parallelepiped, remain
symmetric w.r.t.\ the $xy$-plane. Hence they
remain coplanar.
This symmetry, and consequently coplanarity property remains unchanged if we
still translate vertically both of
these tetrahedra, in a way symmetric w.r.t.\ the
$xy$-plane.

Consider the edges of the already rotated (but not yet
translated)
two tetrahedra, originally lying on one
horizontal face of the rectangular parallelepiped. They
may have
projections to the $xy$-plane, spanning properly
intersecting lines.
In this case, some unique
symmetric vertical translations of the
already rotated tetrahedra will make the lines spanned
by these two already rotated, and translated
edges intersecting. Moreover,
this happens simultaneously for the pairs of edges,
originally lying on both horizontal faces.

Suppose $- \pi \le \varphi \le \pi $.

First we investigate the
case $- \pi < \varphi < \pi $.
We are going to show that in this case we have the following.
The edges of the already rotated (but not yet translated)
two tetrahedra, originally lying on one
horizontal face of the rectangular parallelepiped, satisfy
the following. They have
projections to the $xy$-plane, which span properly
intersecting lines.

The lines spanned by the
diagonals of the upper horizontal face of the rectangular
parallelepiped
divide the plane
spanned by this face to four angular
domains. We translate these angular
domains
to the $xy$-plane, so that their
vertices are translated to ${\bold{0}}$. We may
suppose that ${\bold{u}}$ lies, e.g.,\
in the closed angular domain
$(Q_2 - Q_1){\bold{0}}(P_4 - P_3)$ (in the unrotated
position). We investigate the projections
of the rotated diagonals $Q_1Q_2$ and $P_3P_4$
to the $xy$-plane.
Let us
replace these rotated diagonals by their translated copies
${\bold{0}}(Q_2 - Q_1)$ and ${\bold{0}}(P_4 - P_3)$. Then
the projections of these translated copies
of the rotated diagonals $Q_1Q_2$ and $P_3P_4$,
to the $xy$-plane,
are translates of the projections
of the rotated diagonals $Q_1Q_2$ and $P_3P_4$,
to the $xy$-plane.

Now we investigate the respective rotations of the segments
${\bold{0}}(Q_2 - Q_1)$ and ${\bold{0}}(P_4 - P_3)$,
about the axis
${\bold{0u}}$. Their
points ${\bold{0}}$ remain fixed by these
rotations. Moreover, these segments move on directrices of
two semi-infinite
circular cones, with apex ${\bold{0}}$, and rotation axis 
${\bold{R}}{\bold{u}}$, and with sum of semiapertures less than
$\pi $. (One of these cones,
but not
both, may degenerate to a half-line. This occurs exactly if
${\bold{u}}$  lies on a side of the angular domain 
$(Q_2 - Q_1){\bold{0}}(P_4 - P_3)$, in the unrotated
position.) This rotation axis ${\bold{R}}{\bold{u}}$
divides the
$xy$-plane to two half-planes. Then the projections
of the respective rotated copies of the segments
${\bold{0}}(Q_2 - Q_1)$ and ${\bold{0}}(P_4 - P_3)$
to the $xy$-plane
satisfy the following. They lie,
except for ${\bold{0}}$,
in one open half-plane, and in the complementary closed
half-plane, bounded by this rotation axis, resp.
Hence, by the semiaperture condition,
the lines spanned by them in fact
properly intersect, as stated. Therefore this motion exists
for all $\varphi \in ( - \pi, \pi )$.

There remained to investigate the case $\varphi = \pm \pi $.
Then the edges of the two tetrahedra, originally lying on a
horizontal face of the rectangular parallelepiped, are rotated to positions lying on two
parallel vertical planes. These planes
have as distance the height
$2d_3$ of our rectangular parallelepiped. So no vertical
translations can make these rotated edges of the two tetrahedra
intersecting. However, these vertical translates still
can be coplanar, namely if they are parallel. To check
their possible parallelity,
it is sufficient to replace them
by their already rotated and
translated copies ${\bold{0}}(Q_2 - Q_1)$ and 
${\bold{0}}(P_4 - P_3)$, and check their parallelity. Since
both of them contain ${\bold{0}}$, therefore
their parallelity means
their coincidence. They lie in a vertical plane containing
${\bold{R}}{\bold{u}}$. They coincide exactly if
${\bold{u}}$ lies on the bisector of the angular
domain $(Q_2 - Q_1){\bold{0}}(P_4 - P_3)$, in the
unrotated position. We have  the analogous statements for any
of the above four angular domains.
Since the
horizontal face of the rectangular parallelepiped
is a rectangle, this means that ${\bold{u}}$
(of length $1$)
equals $\pm {\bold{e}}_1$ or $\pm {\bold{e}}_2$. Thus
we have a motion of the first kind, which exists for all
$\varphi $, so, in particular, for $\varphi = \pi $.

For a motion of the intermediate kind (in particular, for a
motion of the first kind, with $\varphi \ne \pi $),
for any given rotation~$\bold A$,
the translation
$\bold b$ is uniquely determined.
This is also shown in the same way as in
\cite{12}, p.~429 (or cf.\ above, at the
``properly intersecting lines'').
However, possibly for such
an $\bold A$ there
exist several other $\bold b$'s yielding a motion of
another kind, see the next paragraph.

We continue with describing the novel motion of the
sixth kind.
Let us suppose that our rectangular parallelepiped satisfies
$d_k = d_i d_j / (d^2_i + d^2_j)^{1/2}$, where $i,
j, k$ are different. 
Here $d_k$ is half the length of the altitude belonging to
the hypotenuse of the right triangle bounded by two sides
and a
diagonal of a face perpendicular to $\bold e_k$.
Let, e.g.,\ $d_3 = d_1 d_2 / (d^2_1 + d^2_2)^{1/2}$.
Let us take an axis of rotation passing through $\bold 0$
and parallel to one
of the diagonals of a horizontal face, say, to $P_1 P_2$,
with axis vector ${\bold{0}}(P_2 - P_1)$.
Let
us consider
a rotation about this axis through the angle
$\pm \pi /2$.

For convenience, rather than taking one tetrahedron as fixed, the other one as
moving, we argue as at the motion of the intermediate kind.
We rotate both tetrahedra about this axis, through an angle
$\pm \pi /4$, in
a way symmetric w.r.t.\ the $xy$-plane. Then we
translate vertically both of
these tetrahedra, in a way symmetric w.r.t.\ the
$xy$-plane, through a suitable distance. Thus we achieve
that all the pairs of the edges
of the tetrahedra, which originally lay on the same
face of the
parallelepiped, will be, simultaneously, coplanar.
Thus we
obtain a position
corresponding to a motion of the intermediate kind.
This is already known to exist, even for all
$\varphi \ne \pi $, so, in particular,
for $\varphi = \pm \pi / 2$.
See Fig.~5,

\noindent
FIGURE 5 ABOUT HERE
%
\smallskip
\vbox{
\centerline{\epsfxsize=81mm 
\epsfbox{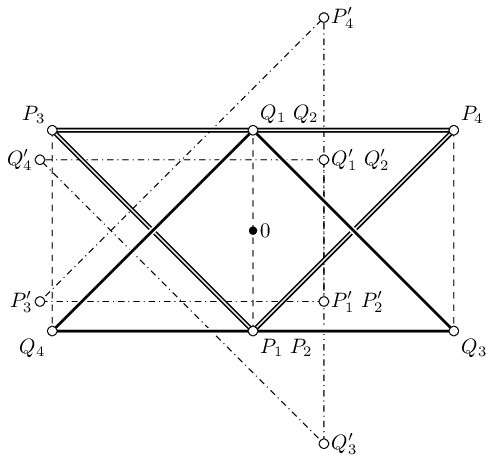}}
\smallskip
\noindent
{
{\centerline{Figure 5.
The motion of the sixth kind.}}
\par}
}
\smallskip
%
\noindent
where the rotated and not yet translated tetrahedra are denoted by
$P_1'P_2'P_3'P_4'$ and
\newline
$Q_1'Q_2'Q_3'Q_4'$. Here $P_i'$ and
$Q_i'$ are the rotated copies of $P_i$ and $Q_i$, resp.
We assume that the sense of rotation is as drawn in Fig.~5 (the other case is
similar, only the role of the indices is changed).

The figure shows the orthogonal projection of the
rectangular parallelepiped along the
axis of rotation.
By $d_3 = d_1 d_2 / (d_1^2 + d_2^2)^{1/2}$
this projection is a rectangle, whose
horizontal side is twice as long as its vertical side.
Note that this implies that the edges $P_1'P_3'$,
$P_2' P_3'$ and their
mirror images $Q_2' Q_4'$, $Q_1' Q_4'$ w.r.t.\ the $xy$-plane
are horizontal. Hence $P_1'P_3'$, $Q_2' Q_4'$, as well as
$P_2' P_3'$, $Q_1' Q_4'$,
are parallel, resp.
Therefore, their any translated copies remain coplanar, resp.
Further the edges $P_1' P_4'$, $P_2' P_4'$, and their mirror images $Q_2'
Q_3'$, $Q_1' Q_3'$ w.r.t.\ the $xy$-plane
lie in a vertical plane (whose projection is a vertical
line in the figure).
Hence any translation in this vertical plane will leave them
coplanar. 
Last, the edges $P_1' P_2'$ (whose projection in the figure is a point) and
$Q_3' Q_4'$ will intersect, thus will be coplanar, after
some
vertical translations of the two tetrahedra, in a way
symmetric w.r.t.\ the $xy$-plane.
Then, by symmetry, also $Q_2' Q_1'$ and
$P_4' P_3'$ will intersect.
However, then any subsequent translation
of the moving tetrahedron in the direction of $P_1'P_2'$ (=~the direction
of $Q_2' Q_1' =$ the direction of the rotation axis) leaves
both of these last two pairs of edges coplanar.

Summing up: each pair of originally intersecting edges
remains coplanar if we make the following operations.
First we make two
rotations --- symmetric w.r.t.\ the $xy$-plane ---
about the above
described axis of rotation, through angles
$\pm \pi /4$. Second we make two
vertical translations --- symmetric
w.r.t.\ the $xy$-plane ---
making $P_1' P_2'$ and
$Q_3'Q_4'$
intersecting. Third we translate the moving tetrahedron in the
direction of the
rotation axis, through an arbitrary distance.
Hence, this is an affine
$1$-ma\-ni\-fold of solutions for $\bold{b}$. (An
{\it{affine
manifold}} is a translate of a linear subspace, or the empty
set.)
We call this, for any
choice of the permutation $(ijk)$ (for which
$d_k = d_i d_j / (d_i^2 + d_j^2)^{1/2}$ holds) a
{\it{motion of the sixth kind}}. (Observe that then
$d_k < d_i,d_j$, hence if this condition holds, then
it can hold only for one $k$.)

Last, we turn to the motion of the third kind, which we are able to give in
analytic form only. Let $\bold u =
[u_1\ \ u_2 \ \ u_3]^T$ and $u_1^2 + u_2^2 + u_3^2 = 1$.
Let $D_i = d_i^{-2}$ and let $\bold{\boldsymbol\Phi(x) = Ax + b}$ be a motion,
where $\bold A$ is a rotation about the axis $\bold{0u}$ through an
angle~$\varphi$.
Note that $(\bold u, \varphi)$ and $(-\bold u, - \varphi)$ represent the same
rotation.
Suppose $0 < \varphi < 2\pi$, with $\varphi \neq \pi /2, 3\pi /2$, and
let $s := \cot
(\varphi /2)$ $(\neq \pm 1)$.
Then, for $i = 1, 2, 3$, let
$$
u_i^2 := \frac{(s^4 + 3s^2 + 1) - (3s^4 + 7s^2 + 1) D_i/(D_1 + D_2 +
D_3)}{2(s^2 + 1)},
$$
provided
$$
\frac{s^4 + 3s^2 + 1}{3s^4 + 7s^2 + 1} > \frac{\max D_i}{D_1 + D_2 + D_3}
$$
(then for each $i$ we have $u_i^2 > 0$).
This determines the rotation part $\bold A$ of the motion
$\bold{\boldsymbol\Phi(x)}$, and the translation part $\bold b$ is then
uniquely determined among all admitted motions. For given
signs of the $u_i$'s, 
this is a $1$-manifold of solutions. 
(But note that actually $\pm ({\bold{u}}, \varphi )$ give
the same rotation.)

We exclude the case $D_1 = D_2 = D_3$, since it has been settled in
\cite{10} and \cite{12}. First
we restrict our attention to the case $u_1, u_2, u_3 > 0$.
Then,
for $\max D_i / (D_1 + D_2 + D_3) \geq 5/11$, this solution
manifold is 
connected. On the other hand, for $\max D_i / (D_1 + D_2 + D_3) < 5/11$,
this solution
manifold consists of three 
connected components, one for $s < -1$ (with $s$ bounded),
one for $-1
< s < 1$, and one for $s > 1$ (with $s$ bounded).

Suppose $D_i \ge D_j,D_k$ (thus earlier
$\max D_i$ will be
replaced by $D_i$).
Here and later till Theorem 1, $(ijk)$ is a
permutation of $\{ 1,2,3 \} $.
For $D_i / (D_1 + D_2 + D_3) \ge 5/11$ this solution
manifold ends (at the infimum or supremum of $s$)
at two points satisfying $u_i = 0$ 
and $|s| \in (0,1]$ (with $|s| = 1$ exactly for
 $D_i / (D_1 + D_2 + D_3) = 5/11$).
For $D_i / (D_1 + D_2 + D_3) < 5/11$ the component for
$s < -1$ begins
(at the infimum of $s$) and the component for $s > 1$ ends
(at the supremum of $s$)
at two
points satisfying $u_i = 0$ 
and $|s| > 1$. 
Moreover, there are no endpoints at $s = \pm
1$, the manifold components going to infinity there.
For any of the above two cases,
these two endpoints are not considered as lying on this
solution manifold. They correspond to two motions of the
intermediate kind.

If $D_i > D_j, D_k$, then the
above mentioned two endpoints (for any of the above two
cases) do not correspond to
motions of the first kind.
If $D_i = D_j > D_k$, then
the above mentioned two endpoints (for any of the above two
cases)
correspond to motions of the first kind.

This motion of the third kind never passes through the basic
position, unless $D_1
= D_2 = D_3$, when it does (for $s^2 = \infty $).

These formulas in fact describe a motion of our
tetrahedra. This fact, and
the mentioned properties of this motion will follow from the
proof of Theorem~1.

For $\min u_i^2 > 0$ this
$1$-manifold of motions of
the third kind is smooth.
However, unless $D_1 = D_2 = D_3$,
this $1$-manifold, parametrized by $ \varphi $, or
by $s$,
is in general only topologically
a manifold with boundary (i.e., its above two endpoints).
Namely, at its
endpoints it is in general not differentiable.
For $D_i \ge D_j,D_k $, and for $s$ yielding
$u_i = 0$, we have
$du_i/ds = (du_i^2/ds)/(2u_i) = \pm \infty $,
provided that $du_i^2/ds \ne 0$.
Here $du_i^2/ds$ exists, and is finite.
E.g., for
$D_i = D_j = 5/11$ and $D_3 = 1/11$, an easy calculation
shows that for $s$ yielding $u_i = 0$
we have $du_i^2/ds \ne 0$. Therefore in general, for this $s$,
we have $du_i^2/ds \ne 0$.
Hence this $1$-manifold with boundary is, in general, not
differentiable at its
endpoints.

Even for $D_i = D_j > D_k$, in general, for $s$ yielding
$u_i = u_j = 0$, we have 
$du_i^2/ds = du_j^2/ds \ne 0$. 
This is shown by the same example.
Similarly, for $D_i \ge D_j , D_k$,
we can specialize to $D_i/(D_1 + D_2 + D_3) = 5/11$, or to
$D_i/(D_1 + D_2 + D_3) = D_j/(D_1 + D_2 + D_3) = 5/11$.
Then we see
that also in these two special cases,
in general, for $s$ yielding
$u_i = 0$, we have 
$du_i^2/ds \ne 0$. 
So, even in each of these three special cases,
we have in general non-differentiability of this
$1$-manifold with boundary, at its endpoints. 

Now we allow any signs of the $u_i$'s. 
The number of $1$-manifolds
of motions of the third kind is $2^3/2
=4$. Namely the number of sign combinations of the $u_i$'s,
when they
are non-zero, is $2^3$, but $\pm ({\bold{u}}, \varphi )$ give
the same rotation.
Then
for $D_i > D_j,D_k$ these four
$1$-manifolds by
twos have both endpoints in common. However,
for different twos
these both endpoints are disjoint pairs.
These endpoints correspond
to motions of the intermediate kind (but not of the first kind).
For $D_i = D_j > D_k$ all four of
these $1$-manifolds have both endpoints in common.
These endpoints correspond to motions of the first kind.

At both endpoints of this motion (at the infimum or
supremum of $s$)
we have bifurcations. First suppose $D_i > D_j,D_k$. Then at
an endpoint the solution set locally consists of two
$1$-manifolds with boundary (this endpoint),
of motions of the third kind.
Moreover, of
one $2$-manifold, of motions of the intermediate
kind.
Second suppose $D_i = D_j > D_k$. Then at
an endpoint the solution set locally consists of all four
$1$-manifolds with boundary (this endpoint),
of the third kind.
Moreover, of one $1$-manifold
of the motions of the first kind,
and of two $2$-manifolds of the motions of the
intermediate kind.

Observe that the motions of the fifth and sixth kinds
exist only
for $\varphi = \pi $, i.e., $s = 0$, and for $\varphi = 
\pi /2, 3 \pi /2$, i.e., $s = \pm 1$, resp.
However, $s = 0$ does not yield an
endpoint of the $1$-manifold of motions of
the third kind. Hence
the motion of the fifth kind cannot occur locally at 
the endpoints of the $1$-manifold of motions of
the third kind.
Further, $s = \pm 1$ yield
endpoints of the $1$-manifold of motions of
the third kind if and only if $\max D_i / (D_1 + D_2 + D_3)
= 5/11$.
The motion of the sixth kind exists only for $D_i =
D_j + D_k$ (where $(ijk)$ is a permutation of
$\{ 1,2,3 \} $; cf.\ {\bf{(7)}} of the proof of Theorem 1).
If an endpoint of the $1$-manifold of
motions of the third kind, and a motion of the
sixth kind coincided, then we had
$1/2 = D_i / (D_i + D_j + D_k) = 5/11$, a contradiction.
Hence
the motion of the sixth kind 
cannot occur locally at 
the endpoints of the $1$-manifold of the motions of
the third kind.

Also at
each {\it{counterbasic position}}, i.e., a motion of the first
kind with rotation angle $\varphi = \pi $ (i.e., $s = 0$),
we have bifurcations. We may suppose ${\bold{u}} = (0,0,1)$.
At a counterbasic position there
locally occur one $1$-manifold of the motions of the
first kind, two $2$-manifolds of the motions of the
intermediate kind, and
one $2$-manifold of the motions of the fifth kind.
However, manifolds
of the motions of the third kind cannot occur at a
counterbasic position. Namely, for
$s = 0$ the formula for $u_i^2$
(at the introduction of the motion of the third kind)
simplifies to $u_i^2 =
\left( 1 - D_i / (D_1 + D_2 + D_3) \right) /2$, for
$i = 1, 2, 3$. Then
$0 = u_1^2 = \left( 1 - D_1 / (D_1 + D_2 + D_3) \right) /2
> 0$, a contradiction. Also
manifolds of the motions of the
sixth kind cannot occur at a counterbasic position,
because for them we have $\varphi = \pm \pi / 2$.

\medskip

We recall that we mean by a motion an isometry (congruence) of the
space, of determinant $+1$.


\medskip

{\bf{Theorem 1.}}
{\it{Consider the two tetrahedra $P_1 P_2 P_3 P_4$ and
$Q_1 Q_2 Q_3 Q_4$, derived
above from the rectangular parallelepiped of vertices $(\pm d_1, \pm d_2, \pm
d_3)$.
The only finite motions admitted by these tetrahedra --- i.e.,
all positions
of the moving tetrahedron, satisfying~{\rm (A)} --- are
the following: The motions
of the first,
intermediate, third, fifth kinds and, provided
$d_k = d_i d_j / (d_i^2 + d_j^2)^{1/2}$ for some permutation
$(ijk)$ of $\{ 1,2,3 \} $, of the sixth kind,
described above.}}

\medskip


{\it{Proof.}}
{\bf{(1)}}
Let $\bold{\boldsymbol\Phi(x) = Ax + b}$ be a finite motion admitted by our
tetrahedra. I.e., $\bold A = [a_{ij}]$ is an orthogonal $3 \times 3$ matrix
of determinant $+1$, and
$\bold b = [b_1 \quad b_2 \quad b_3]^T$ is a vector in
${\bold R}^3$, and condition (A) is satisfied. Here
$\bold A$ is a rotation about some axis $\bold{0u}$, where
$\bold u = [u_1 \quad
u_2\quad u_3]^T$ and $u_1^2 + u_2^2 + u_3^2 = 1$, through an angle
$\varphi$ (with sense of rotation as described in~\S1).

For the vertices $Q_i(x_i,y_i,z_i)$ of the moving
tetrahedron
we have, using the coordinates of the basic position of
the $Q_i$'s from the
beginning of {\bf 2.1},
$$
\align
Q_1(x_1, y_1, z_1) &= \bold A [-d_1 \quad d_2 \quad d_3]^T + \bold b,\\
Q_2(x_2, y_2, z_2) &= \bold A [d_1 \quad {-}d_2 \quad d_3]^T + \bold b,\\
Q_3(x_3, y_3, z_3) &= \bold A [d_1 \quad d_2 \quad {-}d_3]^T + \bold b,\\
Q_4(x_4, y_4, z_4) &= \bold A [{-}d_1 \quad {-}d_2 \quad {-}d_3]^T + \bold b.
\endalign
$$
Let us denote
$$
\bold D = \bmatrix d_1 & 0 & 0 \\
0 & d_2 & 0 \\
0 & 0 & d_3\endbmatrix \!.
$$
Using the notations $P_i^0$, $Q_i^0$ from \S1
(e.g.,\ Fig.\ 3),
and $P_i$ from the beginning of
{\bf 2.1}, we have $P_i = \bold D P_i^0$.
Further we have $Q_i = {\boldsymbol \Phi } (\ol Q_i)$, where $\ol Q_i$, or $\ol
Q_i^0$, is the basic position of $Q_i$, or $Q_i^0$, resp.
These satisfy $\ol Q_i = \bold D \ol Q_i^0$.
The coplanarity, e.g.,\ of the fixed vertices $P_1 = \bold D P_1^0$ and $P_2
= \bold D P_2^0$, and of the moving vertices $Q_3 =
{\boldsymbol \Phi }(\ol Q_3)
= {\boldsymbol \Phi }\bigl({\bold D}(\ol Q_3^0)\bigr)$ and
$Q_4 = {\boldsymbol \Phi }(\ol Q_4) = {\boldsymbol \Phi }
\bigl({\bold D}(\ol Q_4^0)\bigr)$,
is equivalent to the following. The points $P_1^0$, $P_2^0$, and ${\bold
D}^{-1} \bigl({\boldsymbol \Phi }({\bold D}(\ol Q_3^0))\bigr) = {\bold D}^{-1} {\bold
A}{\bold D}(\ol Q_3^0) + \bold D^{-1} \bold b$ and ${\bold
D}^{-1} \bigl({\boldsymbol \Phi }({\bold D}(\ol Q_4^0))\bigr) = {\bold D}^{-1} {\bold
A} {\bold D}(\ol Q_4^0) + \bold D^{-1} \bold b$ are coplanar.
We have analogous
equivalent conditions for the coplanarity of the
other quadruples
of vertices to be considered.

So $\bold{\boldsymbol\Phi(x) = Ax + b}$
represents a motion of our tetrahedra if and only if the
following holds. The
transformation $\bold
D^{-1} \bold{A D(x) + D}^{-1} \bold b$ of the vertices $Q_1^0, \dots, Q_4^0$
preserves coplanarity of the four vertices of any face of the cube with
vertices $(\pm 1, \pm 1, \pm 1)$. Simultaneously,
two of these vertices,
namely the vertices $P_i^0$, are
fixed, and only two of
them, namely the vertices $Q_i^0$, are transformed by this
transformation.
We write
$[a_{ij}^0] := \bold D^{-1} \bold A \bold D =
[d_i^{-1} a_{ij} d_j]$, and
$[b_1^0 \quad b_2^0 \quad b_3^0]^T := \bold D^{-1} \bold b
= [d_1^{-1} b_1 \quad d_2^{-1}b_2 \quad d_3^{-1}b_3]^T$.
Then,
like in \cite{12}, p.~435, by these coplanarities we have
$$
\align
-(a_{22}^0 + a_{33}^0) b_1^0 + a_{12}^0 b_2^0 + a_{13}^0 b_3^0
&= a_{21}^0 a_{13}^0 + a_{31}^0 a_{12}^0 + (a_{23}^0 + a_{32}^0)(1 - a_{11}^0),
\tag{I/1}\\
-(a_{23}^0 + a_{32}^0) b_1^0 + a_{13}^0 b_2^0 + a_{12}^0 b_3^0
&= a_{12}^0 a_{21}^0 + a_{31}^0 a_{13}^0 + (a_{22}^0 + a_{33}^0)(1 - a_{11}^0).
\tag{II/1}
\endalign
$$
There hold the analogous equations obtained from these ones by the permutation
of the indices $1 \mapsto 2 \mapsto 3 \mapsto 1$. These equations will be denoted by (I/2) and
(II/2), resp.
Similarly, using the permutation $1 \mapsto 3 \mapsto 2 \mapsto
1$ we get equations
(I/3) and (II/3).
Thus we obtain a system of six linear equations for $b_1$,
$b_2$ and $b_3$. This system of equations
expresses the coplanarity of the respective six
quadruples from the fixed
points $P_i$ and the moved points $Q_i$, cf.\ \cite{12}, p.
435.

\smallskip

{\bf{(2)}}
In these six equations (I$/i$), (II$/i$), where $i = 1, 2, 3$,
we replace $a_{ij}^0$ by $d_i^{-1} a_{ij} d_j$, and $b_i^0$
by $d_i^{-1}b_i$.
Then we express $a_{ij}$ by $u_1, u_2, u_3$ and $\varphi$, by the well-known
formula
${\bold{Ax}} = \langle {\bold{x}},{\bold{u}} \rangle
{\bold{u}} + \cos \varphi \cdot ( {\bold{x}} -
\langle {\bold{x}},{\bold{u}} \rangle {\bold{u}} )
+ \sin \varphi \cdot ( {\bold{x}} \times {\bold{u}} )$
(cf.,\ e.g.,\ \cite{12}, p.~436).

\smallskip
{\narrower{\narrower\noindent
Then
let us consider the $6 \times 4$ matrix (B$'$) formed by
the coefficients of
$b_1, b_2, b_3$, and the right-hand sides of these
six equations, rewritten as indicated
above.
(Its entries depend on ${\bold{u}}$, $\varphi $, and the constants
$d_i$.)
\par}}


\smallskip
\noindent
We call its rows I/$i$ and II/$i$, according to the equation
they correspond to. We will not write out (B$'$) explicitly,
but rather we will make some simplifications in it.
In matrix (B$'$) we multiply rows I/$i$, II/$i$ by $d_i$, then
divide the fourth column by $2
d_1 d_2 d_3$, and then divide row II/$i$ by
$d_1 d_2 d_3 d_i^{-1}$.
Like in {\bf{2.1}}, before Theorem~1, denote $D_i = d_i^{-2}$.
Note that for $\varphi = 0$, thus $\bold{A = I}$, we have
the unique motion
with ${\bold{b}} = {\bold{0}}$.
Henceforward we will assume $0 < \varphi < 2\pi$.
Letting $s := \cot (\varphi /2)$, as in {\bf{2.1}},
before Theorem 1, we still multiply each
entry of the matrix obtained last time from (B$'$)
by $(s^2 + 1)/2$.
Thus we obtain the matrix 
$$
\bmatrix
u_1^2 - s^2 & u_1 u_2 + su_3 & u_3 u_1 - su_2 & u_2 u_3(D_2 + D_3)\\
u_1u_2 - su_3 & u_2^2 - s^2 & u_2 u_3 + su_1 & u_3 u_1(D_3 + D_1)\\
u_3u_1 + su_2 & u_2 u_3 - su_1 & u_3^2 - s^2 & u_1 u_2(D_1 + D_2)\\
-u_2 u_3 (D_2 + D_3) +  & (u_3 u_1 - su_2)D_2 & (u_1 u_2 + su_3)D_3 & 0\\
+ su_1(D_3 - D_2) & & & \\
(u_2 u_3 + su_1)D_1 & -u_3 u_1(D_3 + D_1) + & (u_1 u_2 - su_3)D_3 & 0\\
& + su_2(D_1 - D_3) & & \\
(u_2 u_3 - su_1) D_1 & (u_3 u_1 + su_2) D_2 & - u_1 u_2(D_1 + D_2) + & 0 \\
& & + su_3(D_2 - D_1) &
\endbmatrix \!.
\tag{B}
$$
(Note that for $D_1 = D_2 = D_3 = 1$ this reduces to (B)
in \cite{12}, p.~437, up to a 
factor $1/2$ in the fourth
column.)
The rows of matrix (B) corresponding to equations (I/$i$), (II/$i$) will
be called rows I/$i$, II/$i$ of~(B). Later, unless stated
otherwise, we will consider
only rows I/$i$, II/$i$ of matrix (B), and not of the
original matrix (B$'$). Now
suppose that ${\bold{u}}$ and $s$ are
fixed.
Then the solvability of the system of equations corresponding
to
this new matrix (B), for ${\bold{b}}$, is
equivalent to the solvability
of the system of equations corresponding to the original
matrix (B$'$) (thus of our original system
of equations (I/$i$), (II/$i$), for
$i = 1, 2, 3$), for ${\bold{b}}$.
More exactly, ${\bold{u}}$, $s$,
${\bold{b}}$ is a solution
of the system of equations corresponding to (B$'$) if and
only if ${\bold{u}}$, $s$, ${\bold{b}}/(2d_1d_2d_3)$
is a solution
of the system of equations corresponding to (B).
Hence, for ${\bold{u}}$ and $s$ fixed, the dimensions
of the solution manifolds, for ${\bold{b}}$, of the two 
systems of equations (if they are not empty) are the
same. Moreover, they can be obtained from each other by
multiplication with a non-zero constant. Later we will not be
interested in formulas for ${\bold{b}}$, therefore,
{\it{unless stated otherwise, we
will use the system of equations corresponding to}} (B).

\smallskip
{\bf{(3)}}
The upper left $3 \times 3$ submatrix of (B) is independent of $D_i$, hence
it is singular in the same case when
it is singular for $D_1 = D_2 = D_3 = 1$.
Observe that the left-hand side of equation (I/1) equals $d^{-1}_1[-(a_{22} +
a_{33})b_1 + a_{12} b_2 + a_{13} b_3]$, and similarly for (I/2), (I/3).
Hence the determinant of the considered $3 \times 3$
submatrix of (B) is a non-zero number
times the determinant $\bigl|a_{ij} - \delta_{ij}(a_{11} + a_{22} +
a_{33})\bigr|$.
By \cite{11}, p.~270 or \cite{12}, p.~438,
this determinant is $0$ if and only if $\varphi = \pm \pi /2$ or $\varphi
= \pi$.
(We do not distinguish between angles differing by multiples of $2\pi$.)
Hence our equations (for matrix (B))
can have a non-unique solution for $\bold b$ only for
$\varphi = \pm \pi /2$ and $\varphi = \pi$, i.e., for $s = \pm 1$ and $s =
0$.

The last three rows of matrix (B) are linearly dependent.
Namely, multiplying row II/$i$ by $D_i$,
and then summing them, we obtain the zero row.

Multiplying row I/$i$ of (B) by $u_i$, and then
summing them, we obtain
$$
\bigl[ u_1(1 - s^2) \quad u_2(1 - s^2) \quad u_3(1 - s^2)\quad u_1 u_2
u_32(D_1 + D_2 + D_3)\bigr].
$$
The corresponding equation
implies that for $s^2 = 1$, i.e., for
$\varphi = \pm \pi /2$, for any solution
of our equations
(for matrix (B)) we have $u_1 u_2 u_3 = 0$.

\smallskip
{\bf{(4)}}
The determinant of the submatrix of (B) formed by its
rows I/1, I/2, I/3,
II/1 is a homogeneous eighth degree polynomial of
$u_1, u_2, u_3$ and $s$, with coefficients polynomials of the
$D_i$'s.
A straightforward but somewhat lengthy calculation
gives that it equals
$$
\align
&\hskip-8mm
su_1 u_2 u_3(D_1 + D_2 + D_3) \{(u_1^2 + u_2^2 + u_3^2) \times \\
&\hskip-8mm
\times [u_1^2(D_3 - D_2) + u_2^2 (-D_2 - D_3) + u_3^2 (D_2 + D_3)] +\\
&\hskip-8mm
+\! s^2 [ u_1^2(3D_3\! -\! 3D_2) \! + \! u_2^2 (-3D_2\! +\! D_3)\! +\!
u_3^2(-D_2\! +\! 3D_3)]
\!+\! s^4(D_3\! -\! D_2) \! \},
\tag{C}
\endalign
$$
{\it{which equals $0$}}.
Now suppose $s u_1 u_2 u_3 \neq 0$. Then the factor
of (C) in braces
is~$0$. Moreover, two more
analogous expressions are equal to~$0$, which are obtained from this
expression by cyclic permutations of the indices. (These arise analogously from
the determinants of the submatrices formed by
rows I/1, I/2, I/3, II/2,
and I/1, I/2, I/3, II/3 of (B), resp.)
These three equations are homogeneous linear in $D_2, D_3$, in $D_3, D_1$, and
in $D_1, D_2$, resp., and can be written as
$$
\align
&\frac{(u_1^2 + u_2^2 + u_3^2)(-u_1^2 + u_2^2 + u_3^2) + s^2(u_1^2 + 3u_2^2 +
3 u_3^2) + s^4}{D_1} =\\
&= \frac{(u_1^2 + u_2^2 + u_3^2)(u_1^2 - u_2^2 + u_3^2) + s^2(3u_1^2 + u_2^2 +
3 u_3^2) + s^4}{D_2} =\\
&= \frac{(u_1^2 + u_2^2 + u_3^2)(u_1^2 + u_2^2 - u_3^2) + s^2(3u_1^2 + 3u_2^2
+ u_3^2) + s^4}{D_3}
\tag{D}
\endalign
$$
(thus they are actually only two equations).

\smallskip
{\bf{(5)}}
First we discuss the case when the first factor of (C), i.e., $s$,
equals~$0$.
Consider the $4 \times 4$ matrix formed by rows I/1, I/2, II/1, and the
sum of $(D_1 + D_2)$ times row II/2 and $D_3$ times row II/3 of our
matrix~(B).
Its determinant is
$$
-u_1 u_2 u_3^2 D_1 D_3(D_1 + D_2 + D_3) \bigl[u_2^2(D_2 + D_3) - u_1^2(D_3 +
 D_1)\bigr] (u_1^2 + u_2^2 + u_3^2),
$$
{\it{which equals}} $0$.
By cyclic permutation of rows I/$i$ and II/$i$ we get two more
similar
equations. Namely the
expressions, obtained from the last expression by cyclic
permutations of the
indices, are equal to~$0$.

These three equations together imply $u_1 u_2
u_3 = 0$, or $u_1^2 : u_2^2 : u_3^2 =
(D_2 + D_3) : (D_3 + D_1) : (D_1 + D_2)$.

First let us suppose, e.g., $u_1 = 0$.
Then the equation corresponding to row I/1 becomes $0 = u_2 u_3(D_2 + D_3)$,
thus $u_2 u_3 = 0$.
Let, e.g.,\ $u_1 = u_2 = 0$ (thus $\bold{u} = \pm
\bold{e}_3$ by $u_1^2 + u_2^2 +u_3^2 = 1$). By $s = 0$ we have $\varphi = \pi$.
Then the pairs of the edges of the two tetrahedra, originally lying on some
vertical face of the rectangular
parallelepiped, are
parallel. Thus they remain coplanar after
any translation of the moving tetrahedron.
The pairs of the edges of the two tetrahedra,
originally lying on some horizontal face
of the rectangular
parallelepiped, are properly
intersecting and horizontal.
Thus they remain coplanar after a translation
of the moving tetrahedron,
through a vector
$[b_1\quad b_2\quad b_3]^T$, exactly when $b_3 = 0$.
Therefore the set of solution vectors
$[b_1\quad b_2\quad b_3]^T$ of our equations
(for matrix (B))
is given by all vectors with $b_3 = 0$. This gives
the motion
of the fifth kind.

Second let us suppose $u_1^2 : u_2^2 : u_3^2 = (D_2 + D_3) :(D_3 + D_1) :
(D_1 + D_2)$. Then by $u_1^2 + u_2^2 + u_3^2 = 1$ we have $u_1 u_2 u_3 \neq
0$.
Observe that in this case the double equality (D) is
satisfied (for $s = 0$).
We will investigate this case further together with the investigation of~(D),
in {\bf{(10)}} (where both cases $s = 0$ and $s \neq 0$ will be allowed). Now we only show
that in this case there is a unique solution $\bold b$
of our 
equations
(for matrix (B)).
(The existence of the solution for $\bold b$
will be showed in {\bf{(10)}}.)
For this consider the $4 \times 4$ matrix from the beginning of {\bf{(5)}}.
Then take
its $3 \times 3$ submatrix consisting of the first three entries of its
first, third and fourth rows.
Its determinant is a homogeneous polynomial of third degree
of the $D_i$'s, with coefficients homogeneous sixth degree
polynomials of the $u_i$'s. To
test if it is zero or not, it suffices to substitute
$D_2 + D_3 := u_1^2$, and $D_3
+ D_1 := u_2^2$, and $D_1 + D_2 := u_3^2$.
Thus this determinant becomes
$$
u_1^2 u_2 u_3(u_1^2 + u_2^2 + u_3^2)^2 (u_1^2 + u_2^2 - u_3^2)(-u_1^2
+ u_2^2 + u_3^2)/8.
$$
Here $u_1^2 u_2 u_3 \neq 0$, as shown above, and each other factor is positive.
E.g.,\ $(u_1^2 + u_2^2 - u_3^3) / (u_1^2 + u_2^2 + u_3^2) = D_3 / (D_1 + D_2
+ D_3) > 0$.
Hence the considered determinant is non-zero, showing
unicity of the solution
for~$\bold b$ (for matrix (B)).

\smallskip
{\bf{(6)}}
Second we discuss the case when the factor $u_1 u_2 u_3$ of
(C)
equals~$0$.
Let, e.g.,\ $u_1 = 0$.
In this case we have the motions of the first,
intermediate, fifth and sixth kinds.
We have to show for $u_2u_3 \neq 0$, that if the motion of
the intermediate kind
does not exist --- that is $\varphi = \pi$, i.e., $s = 0$ ---
then we do not
have any solution.
(Recall the following. The motion of the
intermediate kind exists only for
$\varphi \neq \pi$. The motion of the fifth kind ---
for $u_1 = 0$ --- exists only
for $u_2 u_3 = 0$. The motion of the sixth kind exists only
for $\varphi = \pm \pi /2$.)
However, in {\bf{(5)}}
it has been shown that $s = 0$ and $u_1 = 0$ imply $u_2 u_3 = 0$.

There remains the question of the
unicity of the translation part $\bold b$ of the
motion.
However, in {\bf{(3)}}
it has been shown that our equations (for matrix (B))
can have a non-unique
solution for $\bold b$ only for $s = 0$ and $s = \pm 1$.
The case $s = 0$ has been dealt with in {\bf{(5)}}, and
the question of unicity has been
completely settled there. We will deal with $s = \pm 1$
in~{\bf{(7)}}.

\smallskip

{\bf{(7)}}
We turn to discuss unicity of $\bold b$ for $s = \pm 1$, i.e., $\varphi
= \pm \pi /2$.
Replacing ${\bold{u}}$ by $-{\bold{u}}$ if
necessary, we may assume $s = 1$.
In {\bf{(3)}} it has been shown that $s^2 = 1$ implies $u_1 u_2 u_3 = 0$.
Let, e.g.,\ $u_3 = 0$.
Then our matrix (B) becomes a function only of $u_1$, $u_2$ and $D_1$, $D_2$,
$D_3$.
Loosing homogeneity, we will use $u_1^2 + u_2^2 = 1$.
Thus we see that rows I/1 and I/2 are proportional, they are $-u_2$, and $u_1$
times $[u_2\quad -u_1\quad 1\quad 0]$, resp.
Since $u_1$, $u_2$ are not both zero,
we may replace rows I/1 and I/2 by one row
$[u_2\quad -u_1\quad 1\quad 0]$. Row I/3 is retained.
As mentioned in {\bf{(3)}}, rows II/1, II/2 and
II/3 are linearly dependent, with
non-zero coefficients. Hence we may omit from among them
row~II/3.
Thus we obtain a $4 \times 4$ matrix ${\bold{M}}$.
The question of the dimension of the
solution manifold, for $\bold b$, of the 
equations corresponding to ${\bold{M}}$, at the
considered rotation part of the motion, is equivalent to
the same question
regarding matrix~(B).
By $u_3 = 0$ and $\varphi = \pm \pi /2 \neq \pi$ one
solution always
exists, namely a motion of the intermediate kind.
Thus the dimension of this solution manifold is $3 - r$, where $r$ is the rank
of the matrix
$$
{\bold{N}} = \bmatrix
u_2 & - u_1 & 1\\
u_2 & - u_1 & -1\\
u_1(D_3 - D_2) & - u_2D_2 & u_1 u_2 D_3\\
u_1 D_1 & u_2(D_1 - D_3) & u_1 u_2 D_3
\endbmatrix\!,
$$
obtained by omitting the last column from 
matrix ${\bold{M}}$.

Subtracting the second row of ${\bold{N}}$
from the first one, the first row becomes $[0\quad
0\quad 2]$. Hence $r = 1 + r'$, where $r'$ is the rank of the $3 \times 2$
submatrix ${\bold{N}}'$ of ${\bold{N}}$ at the lower left corner.
If $u_1$ or $u_2$ is $0$, then
we have $r' = 2$, thus $r = 3$. Then we have a unique
solution for~$\bold b$ (for matrix ${\bold{M}}$).
Now let $u_1 u_2 \neq 0$.
The determinants of the $2 \times 2$ submatrices of ${\bold{N}}'$,
obtained by omitting
its first, second or third row, resp., are the following:
$u_1 u_2 D_3(D_1 + D_2 -
D_3)$, and $u_2^2(D_1 - D_3) + u_1^2 D_1$, and $-u_2^2 D_2 + u_1^2 (D_3 - D_2)$.
If any of these expressions is not~$0$, then
we have $r' = 2$, thus $r = 3$.
Then again there is a unique solution of our equations
for~$\bold b$ (for matrix ${\bold{M}}$).
If all these above expressions are equal to~$0$, then
we have (equivalently) $D_3 =
D_1 + D_2$ and
$u_1^2 D_1 = u_2^2 D_2$. Hence, by $u_3 = 0$, we have
$$
u_1^2 = D_2 / (D_1 + D_2) = d_1^2/(d_1^2 + d_2^2)
{\text{ and }} u_2^2
= D_1/(D_1 + D_2) = d_2^2/(d_1^2 + d_2^2).
$$
In this case $r' = 1$, thus $r = 2$, and then the dimension
of the affine
manifold of solutions,
for $\bold b$ (for matrix ${\bold{M}}$, or matrix (B)),
is $1$ (for this rotation part of the motion). 

It remains to show that geometrically this is the motion
of the sixth kind.
Since $D_i = d_i^{-2}$, therefore $D_3 = D_1 + D_2$ means
$d_3 = d_1 d_2
/ (d_1^2 + d_2^2)^{1/2}$.
We have from above $[u_1\quad u_2\quad u_3]^T = [ \pm d_1 /
(d_1^2 + d_2^2)^{1/2} \quad \pm d_2 /
(d_1^2 + d_2^2)^{1/2} \quad 0]^T$ (the $\pm $ signs being
independent). Hence the axis
of rotation of the rotation part $\bold A$ of the motion is parallel to a
diagonal of a horizontal face of our rectangular parallelepiped in its basic
position.
Further, the angle of rotation is $\pm \pi /2$.
This is just the rotation part of the motion of the sixth
kind
(for $k = 3$, cf.\ the introduction of the motion of the sixth
kind, in {\bf{2.1}}, before Theorem 1).
At describing the motion of the sixth kind, we have
exhibited an affine
$1$-manifold $A_1$ of solutions for $\bold b$, with the above
$\bold A$ (for matrix (B$'$)).
This is therefore
a subset of the entire solution manifold $A_2$, for
$\bold b$, with
this~${\bold{A}}$ (for matrix (B$'$)).
Now we have shown that the entire affine manifold of
solutions $\bold b$,
with this $\bold A$ (for matrix (B)), is
exactly $1$-dimensional.
Recall that the solution manifolds, for ${\bold{b}}$,
for matrices (B$'$) and
(B), are obtained from each other by multiplication with a
non-zero constant, cf.\ {\bf{(2)}}. Hence also
$A_2$ is an affine
$1$-manifold, containing the affine $1$-manifold $A_1$.
Hence 
$A_2 = A_1$, and this is the manifold
of the motions of the sixth kind.

\smallskip
{\bf{(8)}}
Recall that a non-unique solution for $\bold b$ (for matrix
(B)) is possible only for $s = 0$,
$\pm 1$ (cf.\ {\bf{(3)}}). These non-unique solutions
have been discussed in {\bf{(5)}} and
{\bf{(7)}}, resp.

For the existence of solutions we have derived in {\bf{(4)}}
the
equation that (C) equals $0$, and some of its consequences.
The case when the first factor of (C), i.e., $s$, equals~$0$, has been settled
in {\bf{(5)}}, except the case when $u_1 u_2 u_3 \neq 0$ and
(D) is satisfied.
The case when the factor $u_1 u_2 u_3$ of (C) equals~$0$, has been settled in
{\bf{(6)}},
except the case of unicity at $s = \pm 1$. This in turn has been settled
in {\bf{(7)}}.
If $su_1 u_2 u_3 \neq 0$, then
we have derived in {\bf{(4)}} equations~(D).

Therefore all that remains is the following. We have
to solve equations (D), where $s$ can be $0$ or
any non-zero number, and $u_1 u_2 u_3 \neq 0$. Moreover, we
have to verify if they are
solutions of our equations (for
matrix (B)).
Recall that by {\bf{(3)}}
$s^2 = 1$ implies $u_1 u_2 u_3 = 0$, hence we will suppose
$s \neq \pm 1$.

\smallskip
{\bf{(9)}}
Now we show that for $u_1 u_2 u_3 \neq 0$ and $s \neq \pm 1$
any solution
of equations (D) is a solution of our
equations for matrix (B).
This is necessary since equations (D) are only consequences
of the
equations for matrix (B). They have not been gained from
these
equations by
equivalent transformations. Moreover, we do not have a geometrical description
of the motion of the third kind, making its existence evident.

First we show that, for $u_1 u_2 u_3 \neq 0$ and $s \neq 0$,
$\pm 1$,
any solution
of equations (D) is a solution of our equations for
matrix (B).
In fact, equations (D) have been derived for $su_1 u_2 u_3 \neq 0$ from the
equation that expression (C) equals~$0$, and from
two other analogous equations. These
express linear dependence of rows I/1, I/2, I/3 and II/$i$ $(i
= 1,2,3)$ of matrix~(B), resp.
For $s \neq 0$, $\pm 1$ the determinant of the matrix formed by the first
three entries of rows I/1, I/2 and I/3 is not $0$, cf.\ {\bf{(3)}}.
Hence (D) expresses linear dependence of rows II/1, II/2
and II/3 on the linearly
independent rows I/1, I/2 and I/3. (Observe that already their first three
entries form linearly independent row vectors.)
Hence (D) implies that the rank of (B) is at most $3$, thus that the four column
vectors of (B) are linearly dependent.
However, {\it{at this linear dependence the fourth column
vector must have a
non-zero coefficient}}. Namely the first three column vectors are linearly
independent. (Their first three entries already form linearly
independent column vectors.)
{\it{This just means that
our equations have a solution for ${\bold{b}}$
(for matrix {\rm{(B)}})}}.

Now let $u_1 u_2 u_3 \neq 0$ and $s = 0$.
(Recall that $s = \pm 1$ has been excluded, cf.\ {\bf{(8)}}.)
We show that also now any solution of equations
(D) is a
solution of our equations for
matrix (B).
For $s = 0$ (D) gives $(-u_1^2 + u_2^2 + u_3^2)/D_1 = (u_1^2 - u_2^2 + u_3^2)
/ D_2 = (u_1^2 + u_2^2 - u_3^2) / D_3$. Let 
their common value be~$\lambda$, say.
Then $u_i^2 = (\lambda D_j + \lambda D_k)/2$, for $i,j,k$
different (hence $\lambda \neq 0$ by
$u_1^2 + u_2^2 + u_3^2 = 1$), thus 
$u_1^2 : u_2^2 : u_3^2 = (D_2 + D_3) : (D_3
+ D_1) : (D_1 + D_2)$.
Then rows I/1, I/2 and I/3 of matrix (B) are all proportional to $[u_1\quad
u_2\quad u_3\quad 2u_1 u_2 u_3 / \lambda]$. Further, rows
II/1, II/2 and II/3 are
linearly dependent by~{\bf{(3)}}.
Hence (D) implies that the rank of (B) is at most $3$.
Thus after some row manipulations some $3 \times 3$ submatrix, contained in
the first three columns, has a non-zero determinant by the
last paragraph of
{\bf{(5)}} (for $s = 0$). Therefore
we have, like at
the case $s \neq 0$, $\pm 1$, that also in this case
{\it{our equations 
have a solution for ${\bold{b}}$ (for matrix {\rm{(B)}})}}.
Cf.\ the italicized text in the previous paragraph.

\smallskip
{\bf{(10)}}
By the last paragraph of {\bf{(8)}},
it remained to solve equations (D) for $u_1 u_2 u_3 \neq 0$, where $s$ can
be any real number different from $\pm 1$.

Using $u_1^2 + u_2^2 + u_3^2 = 1$, equations (D) become
$$
\frac{(1 - 2u_i^2) + s^2(3 - 2u_i^2) + s^4}{D_i} = \lambda ,
{\text{ where }} \lambda {\text{ is independent of }} i \,\,
(i = 1,2,3).
\tag{E}
$$
Solving this for $u_i^2$, we obtain
$$ 
u_i^2 = \frac12 \left(s^2 + 2 - \frac{\lambda D_i + 1}{s^2 + 1} \right) \!.
\tag{F}
$$
Summing these for $i = 1,2,3$, we obtain
$$
1 = \frac12 \left(3s^2 + 6 - \frac{\lambda(D_1 + D_2 + D_3) + 3}{s^2 +
1}\right) \!,
$$
from which we express $\lambda$ and put it into~(F).
Thus we obtain
$$
u_i^2 = \frac{(s^4 + 3s^2 + 1) - (3s^4 + 7s^2 + 1) D_i/(D_1 + D_2 +
D_3)}{2(s^2 + 1)} ,
\tag{G}
$$
provided of course that all these expressions are non-negative.
Actually, by $u_1 u_2 u_3 $
\newline
$\neq 0$, all these expressions have to be positive.
It is easily seen that these expressions actually satisfy (E)
and have
sum~$1$, thus {\it{we have made equivalent transformations}}.

Using (G), the condition $\min u_i^2 > 0$ is equivalent to $f(s^2) := (s^4 +
3s^2 + 1)/(3s^4 + 7s^2 + 1) > \max D_i / (D_1 + D_2 + D_3)$.
Here $s^2 \mapsto
f(s^2)$ strictly decreases in $[0, \infty)$,
with image $(1/3,1]$.
Hence, except the case $D_1 = D_2 = D_3$, when this inequality is satisfied for
all~$s \in {\bold{R}}$, we have that this inequality is satisfied for $s^2 <
f^{-1} [\max D_i / (D_1 + D_2 + D_3)] < \infty$.
(Thus in this
case this solution set is far from the basic position, which is characterized
by $s^2 = \infty$.)
Here $f^{-1}$, defined on $(1/3, 1]$, and strictly decreasing
there, with image $[0, \infty )$,
is the inverse of~$f$, defined on $[0, \infty)$.
We have $f(1) = 5/11$.

Exclude further the case $D_1 = D_2 = D_3$, which has been
completely settled by \cite{10} and \cite{12}.
We may restrict our 
attention to the case $u_1, u_2, u_3 \ge 0$.
We write $s_0 := [f^{-1} \left( \max D_i /
(D_1 + D_2 + D_3) \right) ]^{1/2}$. By
$f(0) = 1 > \max D_i/(D_1 + D_2 + D_3)$, we have $s_0 > 0$.
By positivity of (G), we have $s \in (-s_0,s_0)$.

First suppose $s_0 < 1$
(i.e., $\max D_i / (D_1 + D_2 + D_3) > 5/11$). 
We claim that $s$ varies in $I := (-s_0, s_0)$.
At $\pm s_0$ the solution
of our equations for ${\bold{b}}$ (for matrix (B)) is
unique,
cf.\ {\bf{(3)}}. For $s \in (-s_0,s_0)$ it
also is continuous in $s$ and
${\bold{u}}$, hence in $s$ (since now, by (G), ${\bold{u}}$ is
continuous in $s$).
Recall from {\bf{(9)}} that 
for $u_1u_2u_3 \ne 0$ and $s \ne \pm 1$ any
solution of equations (D) is a solution for our equations
(for matrix (B)). However, as we have seen 
in (G), also using the equivalence of (G) and (E),
each
$s \in (-s_0,s_0) \not\ni \pm 1$
can occur for a solution for our equations
(for matrix (B)). Moreover, for $s \in (-s_0,s_0)$ we have
$u_i^2 > 0$.
Summing up: in this case we have a
connected manifold of solutions. Moreover, adding to it
its two endpoints, at $s =
\pm s_0$, we obtain a topological $1$-manifold with boundary
(these two points), contained in the solution set.

Second suppose $s_0 > 1$
(i.e., $\max D_i / (D_1 + D_2 + D_3) < 5/11$). 
We claim that $s$ varies in
$J := (-s_0,s_0) \setminus \{ -1, 1 \} $.
From above, we have $s \in (-s_0,s_0)$.
Like in the case $s_0 < 1$, each
$s \in J$ actually can
occur for a solution
of our equations for ${\bold{b}}$ (for matrix (B)).
However, 
for $u_1u_2u_3 \ne 0$, we have that
$s = \pm 1$ cannot occur for a solution, cf.\ {\bf{(3)}}.
We are going to show that $u_1u_2u_3 \ne 0$ continues
to hold also for $s = \pm 1$. 
This will imply that for $0, \pm 1 \ne s \to \pm 1$ the
uniquely existing solutions ${\bold{b}}$ of our equations
(for matrix (B)) tend to infinity (recall
that non-uniqueness for ${\bold{b}}$ can occur only for
$s = 0, \pm 1$, cf.\ {\bf{(3)}}).
(In the contrary case, a standard compactness argument
would give that $s = \pm 1$ could occur for a solution.)
We have $f(s_0^2) = \max D_i / (D_1 + D_2 + D_3)$.
Thus, applying (G) for $s = \pm 1$, we have
$$
\min u_i^2 =
\frac{5/11 - \max D_i / (D_1 + D_2 + D_3)}{4/11} =
\frac{f(1) - f(s_0^2)}{4/11} > 0.
$$
Then however, $u_1^2 u_2^2 u_3^2 \ge (\min u_i^2)^3 > 0$,
as asserted.
Summing up: in this case we have a solution manifold of
three connected components: one for $s \in (-s_0,-1)$, one
for $s \in (-1,1)$, and one for $s \in (1,s_0)$. At $s =
\pm 1$ the manifold components go to infinity.
Like for
$s_0 < 1$, also now, adding to the solution manifold its
two endpoints, at $s = \pm s_0$, we obtain
a topological $1$-manifold with boundary (these two points),
contained in the solution set.

Third suppose $s_0 = 1$ (i.e.,
$\max D_i / (D_1 + D_2 + D_3) = 5/11$). Suppose $D_i \ge
D_j,D_k$.
Then for $s = \pm 1$ the only
possibility of non-uniqueness of the solution for
${\bold{b}}$
for our equations (for matrix (B)) is when $D_i = D_j +
D_k$ (cf.\ {\bf{(7)}}). Then however
$1/2 = D_i / (D_i + D_j + D_k) = 5/11$, a contradiction.
Hence for $s = \pm 1$
the solution for ${\bold{b}}$
for our equations (for matrix (B)) is unique.
However, for $s = \pm 1$ it also exists. Namely, then
$u_i = 0$, and $\varphi = \pm \pi / 2$, which is a motion
of the intermediate kind, which exists even for all $\varphi
\ne \pi $, so, in particular, for $\varphi = \pm \pi / 2$.
Like in the first case
above, $s$ varies in $(-s_0, s_0) = (-1, 1)$. By existence
and unicity of the solution for ${\bold{b}}$
(for matrix (B)),
for $s = \pm 1$, 
we have that the rank of the submatrix of (B), consisting of
its first three columns, for $s = \pm 1$,
is $3$. Otherwise said, for $s = \pm 1$,
some $3 \times 3$ submatrix of matrix (B), contained in its
first three columns, has a non-zero determinant (this
submatrix may
depend on $s = \pm 1$).
Then, for $s = \pm 1$,
the system of the corresponding three linear equations
for $b_1,b_2$ and $b_3$ can be uniquely
solved by Cramer's rule, with
the denominator being non-zero. Then,
in some neighbourhood of an endpoint of this solution
manifold, we have the following.
The 
solution vector ${\bold{b}}$ for this submatrix
depends continuously on the
coefficients of these three linear equations, therefore
on ${\bold{u}}$ and $s$.
Hence, in some neighbourhood of an endpoint of this solution
manifold, {\it{taken in the solution set}},
we have the following. There is a unique solution
for ${\bold{b}}$ (for matrix (B)),
which is furthermore continuous in
${\bold{u}}$ and $s$.
Then
the statements about the endpoints of this solution manifold
follow from this.
Summing up: in this case we have a connected manifold of
solutions. Moreover, adding to it its two endpoints, at
$s = \pm s_0 = \pm 1$, we obtain a topological
$1$-manifold with boundary
(these two points), contained in the solution set.
\hfill
{\bf{QED}}

\medskip


{\bf 2.2.}
In \cite{12}, pp.\ 438--440 a slight generalization of the
question of Tompos's tetrahedra has also been considered. Now we
present the corresponding question for the tetrahedra
$P_1 P_2 P_3
P_4$ and $Q_1 Q_2 Q_3 Q_4$, derived from a rectangular
parallelepiped.
In the physical model of these tetrahedra, we have the
following. The bars (edges) of one tetrahedral
frame (of the fixed tetrahedron $P_1 P_2 P_3 P_4$, say) touch the
corresponding bars (edges) of the other tetrahedral frame (of the moving
tetrahedron $Q_1 Q_2 Q_3 Q_4$) from inside (cf.\ \S 1).
Thus the actual physical constraint is only that each edge $P_i P_j$ lies
``inside $Q_j Q_k$'' (here $(ijkl)$ is
any permutation of $\{ 1,2,3,4 \} $).
This can be defined mathematically as follows
(cf.\ \cite{12}, p.~439).

\smallskip
{\narrower{\narrower\noindent For any permutation $(ijkl)$ of
$\{ 1,2,3,4 \}$,
the signed volume of the tetrahedron $P_i P_j Q_k Q_l$ is either $0$, or has
the opposite sign as that of the tetrahedron
$P_i P_j R_k R_l$. Here $R_k
R_l$ is the translate of the segment $Q_k Q_l$ in the basic position (i.e.,
of the segment $\ol Q_k \ol Q_l = (-P_k)(-P_l)$) satisfying
the following. The midpoint of $R_k
R_l$ is the centre of the rectangular parallelepiped in the basic
position.\par}}

\vskip-12.5pt
\hfill (H)

\smallskip
\noindent
We take (H) as the definition of a {\it{generalized motion of
our moving
tetrahedron $Q_1 Q_2 Q_3 Q_4$}} (while $P_1 P_2 P_3 P_4$ is
fixed). We prove


\medskip

{\bf{Theorem 2.}}
{\it{Consider the pair of tetrahedra, derived from a rectangular parallelepiped,
considered in Theorem~1. For them the generalized admitted
finite motions --- i.e., all positions of the moving
tetrahedron (obtained from its basic position by applying
to it an isometry of the space, of determinant $+1$),
satisfying {\rm{(H)}} ---
are identical with the finite motions admitted by them
(described in Theorem~1).}}

\medskip


\noindent
{\it Proof\/} is analogous to that of
\cite{12}, Theorem~2,
pp.\ 439--440. Details cf.\ there, we only indicate the differences.

Let $\bold{\boldsymbol\Phi x = Ax + b}$ be a generalized admitted finite motion.
For $\bold{A = I}$ we have $\bold{b = 0}$.
From now on we suppose $0 < \varphi < 2\pi$.
Observe that now the constraints are expressed by six inequalities
(corresponding to equalities (1g)--(1l) in
\cite{12},
p.~423). Namely three expressions
(the left-hand sides of (1g), (1h), (1i))
are non-negative, three expressions
(the left-hand sides of (1j), (1k), (1l)) are
non-positive.
In \cite{12} the moving vertices  $Q_i$ were 
obtained by the motion $\bold{Ax + b}$, from the respective
points $(\pm 1, \pm , \pm 1)$.
Differing from this,
now the images by ${\bold{D}}^{-1}$ of the
moving vertices $Q_i$ are obtained
by the transformation $\bold D^{-1} \bold A \bold D(\bold x) + \bold
D^{-1} \bold b$, from the respective
points $(\pm 1, \pm 1, \pm 1)$, while $P_i^0 = 
(\pm 1, \pm 1, \pm 1)$ are fixed.
Cf.\ the proof of our Theorem~1, {\bf{(1)}}.

Subtract from a non-negative above expression a non-positive
above
expression, corresponding to pairs of edges which were
originally diagonals of
opposite faces of the rectangular parallelepiped. Then
divide this difference by $2$. Then like
in \cite{12}, p.~439, we obtain the following.
Rather than our equalities
(II/$i$) in {\bf{(1)}} of the proof of Theorem~1, we will have inequalities. Namely the left-hand
sides of (II/$i$) are not less than their right-hand sides,
which are equal to~$0$. 
However, by {\bf{(3)}}
of the proof of Theorem~1, a positive linear combination of
rows II/$i$ of matrix (B) is~$0$. Therefore the same holds
for matrix (B$'$).
Hence, like in \cite{12}, in each of the inequalities,
corresponding to equalities (II/1), (II/2) and (II/3),
we have equalities.

Now recall that the left hand side of each of the
inequalities
corresponding to an equality (II$/i$)
was obtained as follows. It was half the difference of
a non-negative and a non-positive number. This half
difference being
$0$ implies that both of these
non-negative and non-positive numbers are $0$. In other words,
in all the six original constraint inequalities we have
equalities. I.e.,
each pair of edges $P_i P_j$, $Q_k Q_l$ (where $(ijkl)$ is any
permutation of
$\{ 1,2,3,4 \} $) is coplanar.
Thus $\boldsymbol \Phi $ is a motion admitted by our pair of
tetrahedra.\hfill
{\bf{QED}}


\medskip
{\bf 2.3.}
Let us start, rather than with a rectangular parallelepiped,
with a general parallelepiped.
All the diagonals of all of its faces constitute the edges of
two congruent tetrahedra. This position of the two
tetrahedra is called their {\it{basic position}}.
We define the admitted motions as in~(A), but deleting the
word ``rectangular''. We have,
with the notations from {\bf{(1)}}
in the proof of Theorem~1, that $P_i
= \bold D P_i^0$ and $\ol Q_i = \bold D \ol Q_i^0$. However,
now $\bold D =
[d_{ij}]$ is a general non-singular matrix.
In what follows, we will show (in {\bf{3.1}})
that in certain cases the analogues of the motions
for the case of the cube or
the rectangular parallelepiped exist.
Further we prove the generalization of Theorem~2 to the case of general
parallelepipeds. Also we will investigate the unicity of the
solutions of our equations
for~$\bold b$ (in {\bf{3.1}}). In the physical model, the bars (edges) of
the fixed tetrahedron touch those of the moving tetrahedron
from inside (as in Fig.\ 1).

Also now we have for $[a_{ij}^0]  :=
\bold D^{-1} \bold A \bold D$ and
$[b_1^0\quad b_2^0\quad b_3^0]^T := \bold D^{-1} \bold b$
equations (I/$i$), (II/$i$),
$i = 1,2,3$.
Evidently the left-hand sides of (II/1), (II/2) and (II/3) have sum~$0$.
Their right-hand sides have sum $2 {\text{Tr}}\,
(\bold D^{-1} \bold A \bold D) - 2m_2
(\bold D^{-1} \bold A\bold D)$. Here, for any $3\times 3$ matrix $\bold B$,
we write $m_2(\bold B)$ for
the sum of the symmetric $2 \times 2$ subdeterminants 
of~$\bold B$.
We have ${\text{Tr}}\,(\bold D^{-1} \bold A \bold D) = {\text{Tr}}\, (\bold A)$.
We also have $m_2(\bold D^{-1} \bold A \bold D) =
m_2(\bold A)$. Namely these last
two numbers are the coefficients of $-\lambda$ in the
characteristic
polynomials of $\bold D^{-1} \bold A \bold D$ and $\bold A$,
resp. However, these polynomials coincide.
Hence the sum 
of the right-hand sides of equations (II/1), (II/2) and (II/3)
is
the same as for the case $\bold D = \bold I$, i.e.\ $0$.
(Cf.\ \cite{12}, p. 436, (II/1)$'$).
Thus the sum of equations (II/1), (II/2) and
(II/3) is the equation $0 = 0$.
Hence, among our six linear equations for $b_1,b_2$ and $b_3$,
there are at most five
linearly independent ones.

Defining also {\it{for the case of general parallelepipeds the
generalized admitted
finite motions}} as in Theorem 2,  but in 
(H) deleting the
word ``rectangular'', we have


\medskip

{\bf{Theorem 3.}}
{\it{Consider the pair of tetrahedra, derived above from a general parallel\-epiped. For them the
generalized admitted finite motions are identical with the finite motions
admitted by them.}}

\medskip


\noindent
{\it Proof}
is the same as for Theorem~2, using that the sum of
the linear equations (II/1), (II/2) and
(II/3), for $b_1,b_2,b_3$, is the equation $0 = 0$.
\hfill
{\bf{QED}}


\medskip

Further we will deal with the pair of tetrahedra, derived
above from a general parallel\-e\-pi\-ped, in {\bf{3.1}}.


\medskip

\medskip

3. THE MOTIONS OF A PAIR OF TETRAHEDRA DERIVED FROM A
GENERAL PARALLELEPIPED, OF A PAIR OF REGULAR PYRAMIDAL
FRAMES, AND OF A PAIR OF REGULAR TETRAHEDRA WITH CIRCULAR
ARC EDGES


\medskip

{\bf{3.1.}}
We continue the investigation
of the two tetrahedra derived from a general parallelepiped.
Like in the beginning of {\bf{2.1}}, in
{\bf{(1)}} of the proof of Theorem~1, and in
{\bf{2.3}}, our
parallelepiped is taken as the image of the cube with
vertices $(\pm 1, \pm 1, \pm 1)$, by the non-singular
matrix ${\bold{D}}$.

Similarly like in {\bf{(3)}} of the proof of Theorem~1, we
have the following.
A non-unique solution of
equations (I/$i$), (II/$i$), $i = 1,2,3$, for~$\bold b$, can
occur only if
${\text{Tr}}\,(\bold D^{-1} \bold A \bold D)$ is
an eigenvalue of $\bold D^{-1} \bold
A\bold D$. Equivalently, ${\text{Tr}}\,(\bold A)$
is an eigenvalue of $\bold A$, i.e., $\varphi
= \pm \pi /2$ or $\varphi = \pi $ (cf.\
\cite{11}, p.\ 270, or \cite{12}, p.\ 438). 

Observe that for $[a_{ij}^0] := {\bold{D}}^{-1}{\bold{A}}
{\bold{D}}$ and $[b_i^0] := {\bold{D}}^{-1}{\bold{b}}$ the
left hand sides of equations (I/1), (I/2) and (I/3) form
the vector $[a_{ij}^0 - \delta _{ij} (a_{11}^0 + a_{22}^0
+ a_{33}^0)] [b_1^0 \quad b_2^0 \quad b_3^0]^T$. Therefore,
like in the first paragraph of {\bf{(3)}}, and in the second
paragraph of
{\bf{(9)}} of the proof of Theorem~1, for
$\varphi \neq \pm \pi /2$, $\pi $, we have the following.
The vanishing of the determinants of the system of
equations (I/1),
(I/2), (I/3), (II/1), and of the system of
equations (I/1), (I/2), (I/3), (II/2),
is also a sufficient condition for the solvability of 
equations (I/$i$), (II/$i$), for $i = 1, 2, 3$, for~$\bold
b$.
(Recall from {\bf{2.3}} the linear dependence among
equations (II/1), (II/2) and (II/3). Namely, their sum
is the equation $0 = 0$.)

Because of the linear
dependence among our equations, it is to be expected that
there
is a $1$-manifold of solutions.
This exists --- and is a motion of the third kind ---
if the parallelepiped has a threefold rotational
symmetry about a
spatial diagonal. Namely, one
tetrahedron remains fixed. Beginning from the basic position,
the other
one is first rotated about this spatial diagonal. Then
it is translated in the direction of this spatial diagonal,
till the coplanarity conditions become, simultaneously,
satisfied (this position is unique). The angle of rotation
can be arbitrary,
except $\pm \pi / 2$.

Now suppose that the mid-plane between two parallel faces
contains ${\bold{0}}$, and is
a plane of symmetry of the
parallelepiped. Let these faces be horizontal.
Then, as follows from the considerations in
{\bf{2.1}}, we have the motions of the
intermediate kind, for any 
$\varphi \ne \pi $, and we have the motions of
the fifth kind. The axis of rotation ${\bold{0u}}$
lies in this plane of
symmetry, i.e., the $xy$-plane
(and beside this it can have an arbitrary direction), and is
perpendicular
to this plane of symmetry, resp.
In this symmetric case, suppose that moreover the projection of the parallelepiped along
a diagonal of one of the mentioned parallel faces is a
rect\-angle of side ratio
$2 : 1$. (The other diagonal of this face having a larger projection than the
altitude belonging to this face.) Then we have a
motion of the sixth kind.

However, in this symmetric case,
as follows from the considerations in {\bf{2.1}}, the
motion of the intermediate kind may exist also for $\varphi =
\pi $. This happens
exactly in the cases when the rotation axis
${\bold{0u}}$ (supposed to lie in the $xy$-plane) is
parallel to an angle bisector of the diagonals
of a horizontal face. Moreover, then there exists also a new
motion, which we call a {\it{motion of the seventh kind}}.
Namely, as in {\bf{2.1}}, we rotate both tetrahedra about
this rotation
axis, in a way symmetric w.r.t.\ the $xy$-plane, through
angles $\pm \pi / 2$.
Then 
we make 
arbitrary vertical translations of the two tetrahedra, in a
way symmetric w.r.t.\ the $xy$-plane. Thus we
get positions satisfying our constraints. For the
case of a rectangular parallelepiped, this yields a
motion of the fourth kind.

For the general case, suppose that the parallelepiped is
nearly a cube --- more exactly, $\bold D$ is
near to $\bold I$.
Let us choose a point of
a solution manifold of the motions
of the third kind for the cube (with ${\bold{D}} =
{\bold{I}}$). Let it correspond to 
${\bold{u}}^0 = (\pm 1/\sqrt{3}, \pm 1/\sqrt{3},
1/\sqrt{3}) $
(the $\pm $ signs being independent), and to a fixed
$\varphi \in (0, 2 \pi )$, with
$\varphi \neq \pi / 2, \pi , 3 \pi / 2$ (with some unique
${\bold{b}}$). 
Then for any fixed $\varepsilon > 0$, for a
sufficiently small perturbation ${\bold{D}}$ of ${\bold{I}}$,
we have the following. 
There is a solution 
for the parallelepiped associated to ${\bold{D}}$,
with $u_1^2 + u_2^2 + u_3^2 = 1$ and
$\| (u_1,u_2,u_3) - {\bold{u}}^0 \| < \varepsilon $,
and with this fixed value of $\varphi $ (with some unique
${\bold{b}}$). 
(Possibly these points $(u_1,u_2,u_3)$,
for $\varphi $ varying in a fixed
closed subinterval of $(0, 2 \pi )$, avoiding
$\pi / 2, \pi , 3 \pi / 2$, and for a
sufficiently small perturbation ${\bold{D}}$ of ${\bold{I}}$,
form a smooth $1$-manifold with boundary.)

In fact, our problem is now equivalent to solving the system
of the two
determinantal equations, mentioned in the third paragraph of
{\bf{3.1}}, for $u_1, u_2, u_3$ and $s$.
For the case of the cube (with ${\bold{D}} = {\bold{I}}$),
these equations say that non-zero multiples of $u_2 -
u_3$, and of $u_3 - u_1$, are $0$, cf.\
\cite{12}, p.~437, 3.
Thus these multiples change their signs at the curves
on $S^2$, given by $u_2 = u_3$ and $u_3 = u_1$, resp.
Hence, after a small perturbation of the equations
(i.e., for ${\bold{D}}$ near to ${\bold{I}}$),
the zero-sets of the perturbed multiples will
be near the above two curves, resp.
Therefore, for each fixed $\varphi \in (0, 2 \pi )$,
where $\varphi \neq \pi /2, \pi , 3 \pi /2$,
and for a sufficiently small perturbation, the following
holds.
We have a solution of this perturbed system of
our two equations,
with $u_1^2 + u_2^2 + u_3^2 = 1$ and
$\| (u_1, u_2, u_3) - {\bold{u}}^0 \| < \varepsilon $, and
with this fixed
value of $\varphi $ (with some unique ${\bold{b}}$).
(Cf.\ \cite{3}, p.~40, Proposition~D.)

Of course, for a general non-singular
${\bold{D}}$,
the basic position is a solution as well.
By the linear dependence of equations (I$/i$) and (II$/i$)
(cf.\ {\bf{2.3}}),
it
is to be expected that the basic position
lies on a solution manifold of dimension at least $1$.
Experiences with models, far from the rectangular
parallelepipeds (more exactly, with all edge lengths rather
different), seem to confirm this, even
with dimension exactly $1$. (Observe that, for a rectangular
perallelepiped, each opposite pair of
edges of the tetrahedra have  equal lengths.)
However, we cannot identify this (assumed)
solution manifold, which we can call the {\it{motion of the
eighth kind}}.

We have determined the infinitesimal degree of freedom, at
the basic position,
for several incongruent parallelepipeds having a threefold
rotational
symmetry about a spatial diagonal. (Details of this
calculation will be given in {\bf{3.2}}.) Except for Tompos's
tetrahedra, this infinitesimal degree of freedom always turned
out to be~$1$.

Hence we have the following.
The infinitesimal degree of freedom, at a point of a
solution manifold, which can be reached from the basic
position by a continuous motion, always satisfying the
constraints, is probably, in general, not greater than~$1$.
Moreover, in general, the basic position does not lie on
a smooth $2$-manifold of solutions.

The simplest unsolved case
is probably that of
a parallelepiped $P$ having a
threefold rotational symmetry about a spatial diagonal. 
Using analogous notations as in Fig.\ 3, let the vertices
be $P_i$ (fixed) and $Q_i$ (moving).
We define
the admissible motion, obtained from
the basic position of the moving tetrahedron
$Q_1Q_2Q_3Q_4$, by a rotation about
the axis of rotation
$P_4Q_4$, say, through an angle $\pm \pi /3$, and a
subsequent
(unique) translation in the direction of this axis. This is
denoted by $\Phi _4^{\pm }$.
Then $\Phi _4^{\pm }$ degenerates $P$ to a regular
double pyramid, with base the triangle
$P_1P_2P_3 = \Phi _4^+(Q_1)\Phi _4^+(Q_2)\Phi _4^+(Q_3)
= \Phi _4^-(Q_1)\Phi _4^-(Q_2)\Phi _4^-(Q_3)$. Moreover,
the plane spanned by this triangle is
a plane of symmetry of this double pyramid. Also, $\Phi _4^+$
and $\Phi _4^-$ lie on a connected component of a
$1$-manifold of solutions, of motions of the
third kind (cf.\ above).

Now let, e.g., $i = 1$. Let $S_1$ be the symmetry w.r.t.\ the
plane spanned by the face $P_2P_3P_4$.
Then the symmetric double pyramid
$(P_1P_2P_3P_4) \cup S_1(P_1P_2P_3P_4)$ 
is a degenerate image of $P$.
It can be obtained by choosing any of the two
orientation-preserving isometries $\Phi _i^+$ and $\Phi _i^-$,
which are also admitted motions, and which are
defined as follows. We have
$\Phi _1^+ (Q_1) := P_2$ and
$\Phi _1^+ (Q_2) := P_3$ and $\Phi _1^+ (Q_3) := S_1(P_1)$
and $\Phi _1^+ (Q_4) := P_4$. Similarly, we have
$\Phi _1^- (Q_1) := P_3$ and
$\Phi _1^- (Q_2) := S_1(P_1)$ and $\Phi _1^- (Q_3) := P_2$
and $\Phi _1^- (Q_4) := P_4$.
Analogously we define $\Phi _i^+$ and $\Phi _i^-$,
for $i = 2,3$.
Then $\Phi _i^+$ and $\Phi _i^-$, for any fixed $i \in \{
1,2,3\} $,
are probably points of assumed three
analogues of the $1$-manifold of
the motions of the third kind, resp.
Moreover, a model
experiment indicates three connected
smooth solution $1$-manifolds, containing
these three pairs of points, resp.

Concluding: it is to be
awaited that the solution manifolds have in general
dimension $1$, and also the infinitesimal degree of
freedom at their points is in general~$1$. Let ${\bold{D}}$
be a
small perturbation of ${\bold{I}}$. Then probably, in some small
neighbourhoods of the four $1$-manifolds of the motions of
the third kind for the cube (with ${\bold{D}} = {\bold{I}}$),
there are four $1$-manifolds of solutions. Possibly these
exist even for each non-singular
${\bold{D}}$ (as they do for the
rectangular parallelepipeds).
Moreover, probably
there is one $1$-manifold of solutions, passing through the
basic position. (For the case of a rectangular
parallelepiped, this may degenerate to have length $0$.)

\medskip


{\bf{3.2.}}
We consider two
congruent right pyramids, with regular $n$-gonal
bases $(n \geq 3)$. Suppose that their axes of rotation
coincide, and the basic edges of one
pyramid intersect the lateral edges of the other one, and
also conversely, 
with the vectors from the centres of
the bases to the respective apices being opposite.
(Without this oppositeness property, the moving pyramid could
coincide with the fixed pyramid.)
Additionally,
we suppose that the direction of some basic edge
of one pyramid
and the direction of some basic
edge of the other pyramid enclose an angle $\pi /n$ (Fig.~6).
This position

\noindent
FIGURE 6 ABOUT HERE
%
\smallskip
\vbox{
\centerline{\epsfxsize=101mm 
\epsfbox{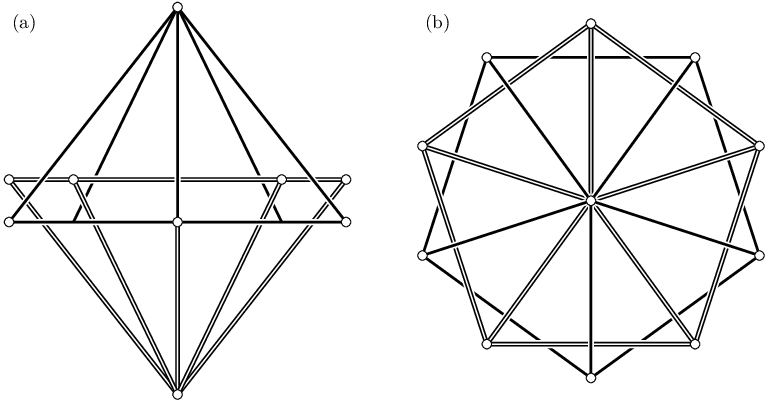}}
\smallskip
\noindent
{
{\centerline{Figure 6.
A pair of regular
pyramids in the basic position:
(a) front view, (b) top view.}}
\par}
}
\smallskip
%

\noindent
is called the {\it basic position\/} of this bar
structure, consisting of these two pyramids.
(Observe that the case $n = 3$ is a special case
of {\bf{3.1}} as well.)

\smallskip
{\narrower{\narrower\noindent
Consider these pyramids as bar structures only. Move each
vertex of the
bar structure consisting of these two pyramids
under the following condition.
One triangular face of one pyramid remains fixed, each bar
(edge) retains its length, and each pair of originally
intersecting 
edges, one from each pyramid, remains coplanar.
\par}}

\vskip-12.5pt
\hfill (I)

\smallskip
The physical model is built in such a way that
the bars (edges) of one
pyramidal frame touch those of the other one from
inside (as in Fig.\ 1).
Observe that the two pyramids are not supposed a
priori to undergo rigid motions (isometries of the space of
determinant~$+1$), but the bases are allowed to change
their shapes, and also to become non-planar. However, by
experimenting with the
respective physical models,
{\it{for positions attainable from the
basic position by continuous motions, always satisfying the
constraints}}, the following seems
probable.
These conditions seem to enforce
the rigid motion of the two pyramids, even with their axes of
rotations coinciding, and with the vectors from the centres
of
the bases to the respective apices being opposite (as in
Fig.\ 6).
Actually even the weakening of (I), 
analogously as in (H), seems to enforce this. Also
cf.\ \cite{9}.

We make local investigations. Since the motions of Tompos's
tetrahedra have already been described, 
we further exclude the case that $n = 3$ and the two
tetrahedra are regular.

\smallskip
{\narrower{\narrower\noindent
We have a
$1$-manifold of finite motions, 
which 
are conjectured to be the only positions, attainable from
the basic position by continuous motions, always satisfying
the constraints.
\par}}


\smallskip\noindent
(This is supported by experimenting with
the models.)
Namely, one pyramid remains
fixed, and the 
other one
undergoes a rigid motion, as follows. 
Its axis of rotation remains fixed, and it undergoes a
certain rotation about this axis, followed by a
suitable translation. This translation happens in
the direction of the common axis of
rotation, through a distance
depending on the angle of rotation. We translate the moving
pyramid
till the coplanarity conditions become, simultaneously,
satisfied (this position is unique).
This is an analogue of the motion of the third
kind for Tompos's tetrahedra. The angle of rotation can
be arbitrary, except $\pm \pi /2$. For $n \ge 4$ (unlike
as for $n = 3$, cf. {\bf{3.1}}), we are
unaware of any other motions,
admitted by our bar structure. 

We have considered
the basic position of this motion, for $3 \leq n \leq 7$,
and for
several different values of
the quotient of the lengths of the lateral and the
basic edges. We have determined for these cases
the infinitesimal degrees of freedom of our 
bar structure, consisting of these two pyramidal frames,
as follows.

The number of the free parameters, i.e., of
all the three
coordinates of all but the fixed three vertices, is $6n - 3$. 
The constraints are that the lengths of all but
the fixed three edges are fixed, and $2n$ pairs of edges
are coplanar, i.e., the tetrahedra spanned by their vertices 
have fixed signed volumes, namely $0$. The total number of
constraints is
also $6n - 3$. Thus we have a function ${\bold{R}}^{6n - 3}
\to {\bold{R}}^{6n - 3}$. This maps a $(6n - 3)$-tuple of
the coordinates of the non-fixed vertices
to the vector with coordinates the $6n - 3$ constraints,
as functions of the previous $6n - 3$ coordinates.
(The constraints are the
lengths of the non-fixed edges, and the signed volumes of
the above tetrahedra.)
Then the number of infinitesimal degrees of freedom
of our bar structure is the
nullity of the Jacobian $J$ of this map, i.e.,
$6n - 3 - {\text{rank}}\,J$.

Having performed these calculations, like e.g.\ in
\cite{7}, \S 2,
we have found the following.
This infinitesimal degree of freedom, at the basic position,
is in all the cases considered by us equal to~$1$,
except in the case of Tompos's tetrahedra. This can be
considered as numerically supporting the above conjecture.

It would be interesting to clarify even that case, when the
two pyramids move as rigid bodies. (Observe that for
$n \ge 4$
this yields an overdetermined system of equations, namely
we have $2n$ equations about coplanarities,
for six unknowns.)


\medskip

{\bf{3.3.}}
Another generalization of the pair of
tetrahedra, derived in \S 1
from the cube with vertices $(\pm 1,
\pm 1, \pm 1)$, is the following.
We replace each edge of both tetrahedra by congruent circular
arcs of some fixed radius. These have 
the same endpoints as the respective edges. Moreover, each
of them lies in the plane spanned by
the respective edge and the centre of the cube.
Further suppose that any congruence of the above
cube to itself
is also a congruence of this system of circular arcs.

Thus we obtain a figure
roughly resembling Fig.\ 1 or Fig.\ 3.
Thus both tetrahedra become tetrahedron-like frames. The
arcs of circles,
replacing diagonals of the same face of the cube,
intersect. Their point of intersection lies on
the straight line,
connecting the midpoint of the cube with the midpoint of the
considered face of the cube. 

\smallskip
{\narrower{\narrower
Fixing one of these frames, we move the other one in the
following way. Each
pair of circles,
containing the pairs of circular arcs, originally
corresponding to
diagonals of some face of the cube, 
continue to have a common point, in the complex
projective sense.
\par}}

\vskip-12.5pt
\hfill (J)

\smallskip
\noindent
Namely, in this sense the
condition is to be awaited simpler.

Now suppose that the frames both lie on the surface of the
circumsphere of the cube. Then an
arbitrary rotation about the centre of the cube, with
translation part
$\bold{b = 0}$, is an admitted motion. So now
we have an at least 3-parameter set of motions.

Again we turn to the general case. We will show that the
motions of the intermediate, third and fifth kinds generically
exist. (These contain the motions of first, second and fourth
kinds as special cases.)

We begin with the analogue of the motion of the fifth kind.
At this motion the moving tetrahedron undergoes from the 
position of first kind --- obtained by ${\bold{A}}$ being
a rotation
through the angle $\pi $, about the
$z$-axis, say --- a translation, through a vector
${\bold{b}}$.
This happens in the following way. An
arbitrarily fixed point of the moving circle,
containing the circular arc corresponding to the edge $Q_1^0 Q_2^0$, in its
rotated position, will coincide after translation with an arbitrarily fixed
point of the fixed
circle, containing the circular arc corresponding to the edge
$P_3^0 P_4^0$ (in analogy with Fig.~4e).

This is a two-parameter motion. At this motion the
circles, containing the
arcs corresponding to the edges $P_1^0 P_2^0$ and
$Q_3^0 Q_4^0$ (in its rotated position), 
also intersect. This follows by a simple
argument using central symmetry. However, the set of
${\bold{b}}$'s
for this ${\bold{A}}$ does not form an affine $2$-manifold ---
on the contrary, it is bounded.
(For
Tompos's tetrahedra we had here an affine $2$-manifold.)
All other pairs of respective circles, which should have
common points, in the complex projective sense, lie in
respectively parallel or coincident planes. This guarantees that these
pairs of circles in fact have common points, in the complex
projective sense. (Observe that the complex projective
extension of any circular line, lying in a horizontal plane,
contains the points $(1, \pm i, 0, 0)$ of the complex
projective space.)

However, for motions of the intermediate and third kinds,
there is a difference as compared to {\bf{2.1}}.

We turn to the analogue of the motion of the third kind.
It will be convenient to rotate our original cube about
${\bold{0}}$, so that $P^0_4$ becomes $(0,0,{\sqrt{3}})$.
The fixed circle $C_{14}$,
containing the circular arc replacing the edge
$P^0_1P^0_4,$ should lie in the vertical plane $y = 0$.
However, the moving circle $C'_{23}$ (obtained by a
rotation,
through an angle $\varphi $, about the $z$-axis,
from its basic position) 
containing the circular arc replacing the edge
$Q^0_2Q^0_3$, will lie in a not vertical plane. 

We denote the projection map to the $xy$-coordinate plane
by $\pi $.
Recalling the definition of the
motion of the third kind, we want to find a $\lambda $,
such that
$C_{14} \cap (C'_{23} + \lambda {\bold{e}}_3) \ne \emptyset
$. (Then, by reason of symmetry, with the evident notations, 
also 
$C_{23} \cap (C'_{14} + \lambda {\bold{e}}_3) \ne \emptyset
$.)

We have that $C_{14}$ lies in the $xz$-coordinate plane,
hence $\pi (C_{14})$ is
contained in the $x$-coordinate axis. The equation system of
$C_{14}$ is $y = 0$, and an equation
of the form $(x-a)^2 + (z - b)^2
= R^2$. Hence, in the complex case, $\pi (C_{14})$ also
contains the 
$x$-coordinate axis, hence equals the 
$x$-coordinate axis.
Even generically $x{\bold{e}}_1$ is
the projection, by $\pi $,
of two points $\left( x, 0, z_1(x) \right),\left(
x, 0, z_2(x)\right) \in C_{14}$ ($C_{14}$
meant in the complex sense): namely the two endpoints of a
chord of $C_{14}$, whose difference lies in the $z$-coordinate
axis.

On the other hand, $C'_{23}$ (meant in the complex sense)
does not lie in a vertical plane, hence the
projection from the affine hull of $C'_{23}$
to the $xy$-coordinate plane is a bijection. Hence
$\pi (C'_{23})$ is a non-degenerate
ellipse, of an equation of the form $Ax^2 + Bxy + Cy^2 + Dx +
Ey + F = 0$. Its intersection with the $x$-axis has an
equation
$Ax^2 + Dx + F = 0$, which has generically two zeroes $x_1,
x_2$. 
Therefore, generically, the number of $(x,0,z)$'s, for which
$(x,0,z) \in C_{14} \cap (C'_{23} + \lambda {\bold{e}}_3)$,
for some
$\lambda $ (depending on $x$ and $z$),
is four. Namely, $x_i{\bold{e}}_1$ is the
projection, by $\pi $, of one point
$\left( x_i,0,z'(x_i)
\right) \in C'_{23}$, and then we choose $\lambda _{ij} =
z_j(x_i) - z'(x_i)$ (for $i,j = 1,2$).
Moreover, these $\lambda _{ij}$'s 
are generically different. This can be seen
from the example with radius of the circles
${\sqrt{2}}$, and $\varphi :=
\pi /2$, where these $\lambda _{ij}$'s are all different.

We turn to the analogue of the motion of the intermediate
kind.
Let the axis of rotation ${\bold{0u}}$ lie in the $xy$-plane. 
Let us apply symmetric rotations, w.r.t.\ the $xy$-plane,
through angles $\pm \varphi /2$, about the rotation axis,
to the two tetrahedron-like frames.
Then the pairs of
the curved edges (circles), corresponding to the two
diagonals
of any of the
originally vertical faces of the cube, remain symmetric
images of each other w.r.t.\ the $xy$-plane.
The (possibly complex projective) intersection points
of the $xy$-plane and one of the circles
lie also on the other circle.

Moreover, after
these symmetric rotations, we have the following.
Each pair of the curved edges,
corresponding to the two diagonals
of an originally horizontal face of the cube, generically
have, as
projections to the $xy$-plane, two elliptical lines. These
projections
have generically four intersection points, in the complex
projective sense. Therefore generically,
for any of the four vertical
lines, containing some of these four
intersection points, the following
holds. Suitable vertical translations,
symmetric w.r.t.\ the $xy$-plane  --- in general through
different 
distances for different lines (cf. below)
--- produce common points of these two curved
edges, lying on this particular vertical line. Further, these
distances, for the two originally horizontal faces, pairwise
coincide. Namely, by reason of symmetry, the distances,
associated to the same vertical line,
(pairwise) coincide, so that these common points, on any of
these, generically four vertical lines
are produced, for the two originally horizontal faces, 
simultaneously.

It remained to show that the four vertical translations
are generically different. Again, it suffices to give one
example for this. We let the radius of the circles to be
${\sqrt{2}}$, and
${\bold{u}} := (1/\sqrt{2}, 
1/\sqrt{2}, 0)$, and $\varphi := \pi $. Then, as at
the analogue of the motion of the third kind, the moving
circle,
containing the arc replacing the diagonal $Q_1^0 Q_3^0$, has
as a projection to the $xy$-plane, in the complex sense, a
fixed
straight line. Therefore we can repeat the considerations
at the analogue of the motion of the third kind, obtaining
that the four $\lambda _{ij}$'s are generically different. 

This question perhaps could be handled in analogy
with
\cite{12}, Theorem 1.
We have six pairs of circles in the space,
one fixed and one moving, which pairwise intersect (in the
complex projective sense).
This for each pair means an
equation of degree six for the coefficients of the
equations of
our circular lines, hence for our parameters $a_{ij}$ and
$b_i$. Thus we have a system of six equations of degree
twelve
for our parameters ${\bold{u}}$, $s$ and
${\bold{b}}$. Unfortunately this is not linear in the
$b_i$'s. Namely, for the analogues of the intermediate
and the third kinds of motions,
there are generically four solutions for ${\bold{b}}$
(in the complex projective sense, cf.\ above). 
Therefore, rather than calculating determinants, as in
\cite{12}, one needs to calculate resultants of polynomials
(cf.\ \cite{14}). Possibly
some symbolic algebraic calculations, like with
Mathematica or Maple,
and efficient algorithms from computational algebraic
geometry could help.


\medskip

{\bf{3.4.}}
In \cite{9} the following general model
was considered, which contains the examples in {\bf{3.1}} and
{\bf{3.2}} as special cases. Let us have two convex polyhedra,
which are combinatorially dual. Let one of them have $f$
faces, $e$ edges and $v$ vertices. Then the other one has
$v$ faces, $e$ edges and $f$ vertices.

\smallskip
{\narrower{\narrower
Consider these polyhedra as bar structures only.
Move each vertex of the bar structure consisting of these two
polyhedra under the following condition. Each bar (edge)
retains its
length, and each pair of combinatorially corresponding
edges is coplanar.

\par}}

\vskip-12.5pt
\hfill (K)

\smallskip
\noindent
The faces may
change their shapes, and may become non-planar, and also
convexity may not hold any more.
The number of free parameters
(all three coordinates of all vertices) is
$3v + 3f$. The number of constraints (edge lengths and
coplanarity conditions) is $3e$. By Euler's theorem,
these numbers have a
difference $6$, i.e., the number of parameters of all rigid
motions of the space.
This used to indicate that there are not
even infinitesimal motions. However, in example {\bf{3.2}},
and sometimes (possibly always) in example
{\bf{3.1}}, there are finite motions, so here intuition
fails.

An example is a pair of congruent tetrahedra with the same
orientation, with the combinatorially
corresponding pairs of edges being those
induced by a fixed orientation-preserving
congruence.
Then a rotation of the
fixed tetrahedron about any of its altitudes,
or a translation of the fixed tetrahedron by any
vector,
yields a moving tetrahedron satisfying the constraints.

\medskip 


{\bf{Problem 1.}} Is there some general theorem behind these
examples, that under suitable hypotheses, the model described
in {\bf{3.4}} always has a finite motion?

\medskip


{\bf{Problem 2.}} Determine the finite motions of the
examples in {\bf{3.1}}, {\bf{3.2}} and {\bf{3.3}}.

\medskip


{\bf{Acknowledgements.}} The authors express their thanks
to Andr\'as Lengyel and Ampar L\'opez for help with the
photography, and special thanks to Andr\'as Lengyel for
producing the line figures of this paper.



\medskip

\medskip

{\centerline{REFERENCES}}

\medskip



[1]
Chen, H.-W.,
{\it{Kinematics and introduction to dynamics of a movable pair of tetrahedra,}}
M. Eng. Thesis, Dept. Mech. Engng., McGill University,
Montreal, Canada,
1991.

[2]
Fuller, R. B.,
{\it{Synergetics. Exploration in the geometry of thinking,}}
Macmillan,
New York,
1975.

[3]
Hurewicz, W., Wallman, H.,
{\it{Dimension theory,}} {\rm Princeton Math. Series,
Vol.~{\bf{4}},}
Princeton Univ. Press,
Princeton, N. J.,
1941.
MR{\bf{3,}}{\rm{312b.}}

[4]
Hyder, A., Zsombor-Murray, P. J.,
{\it{Design, mobility analysis and animation of a double equilateral
tetrahedral mechanism,}} {\rm CIM-89-15 McRCIM Internal Report,}
McGill University,
Montreal, Canada,
1989.

[5]
Hyder, A., Zsombor-Murray, P. J.,
{\it{Design, mobility analysis and animation of a double equilateral
tetrahedral mechanism,}}
{\rm{Proc. Internat. Symp. on Robotics and Manufacturing,}}
ASME Press series, Vol.~{\bf{3}}, 
ISSN 1052--4150,
1990, 49--56.

[6]
Hyder, A., Zsombor-Murray, P. J.,
{\it{An equilateral tetrahedral mechanism,}}
J. Robo\-tics and Autonomous Systems,
{\bf{9}},
(1992),
227--236.

[7]
Kov\'acs, F., Heged\H us, I., Tarnai, T.,
{\it{Movable pairs of regular polyhedra,}}
{\rm{Structural Morphology towards the New Millenium,
Internat. Colloq. Univ. Nottingham, Aug. 15-17, 1997
(Eds. J. C. Chilton, B. S. Choo, W. J. Lewis, O.
Popovi\'c),}} 
Univ. Nottingham, School of Architecture,
Nottingham, UK,
1997, 123-129.

[8]
Makai, E. Jr., Tarnai, T.,
{\it{Overconstrained sliding mechanisms,}}
{\rm{IUTAM-IASS Symp. on Deployable Structures: Theory and
Appl., Proc. IUTAM Symp., Cambridge, UK, 6-9 Sept. 1998 
(Eds. S. Pellegrino, S. D. Guest),}} 
Kluwer,
Dordrecht etc.,
2000, 261-270.

[9]
Pedersen, J., Tarnai, T.,
{\it{Mysterious movable models,}}
Math. Intelligencer
{\bf{34}} {\rm{(3) (2012), 62-66.}}
MR {\bf{2973524}}{\rm{.}}

[10]
Stachel, H.,
{\it{Ein bewegliches Tetraederpaar (A movable pair of tetrahedra, German),}}
Elem. Math.,
{\bf{43}}
(1988),
65--75.
MR {\bf{89i:}}{\rm{51029.}}

[11]
Tarnai, T., Makai, E.,
{\it{Physically inadmissible motions of a movable pair of tetrahedra,}}
{\rm{Proc. Third Internat. Conf. on Engineering Graphics and Descriptive Geometry
(eds. S. M. Slaby and H. Stachel), Vol.~{\bf{2}},}}
Technical Univ.,
Vienna,
1988, 264--271.

[12]
Tarnai, T., Makai, E.,
{\it{A movable pair of tetrahedra,}}
Proc. Royal Soc. London {\rm A}
{\bf{423}},
419--442,
1989.
MR {\bf{90m:}}{\rm{52010.}}

[13]
Tarnai, T., Makai, E.,
{\it{Kinematical indeterminacy of a pair of tetrahedral
frames,}}
Acta Techn. Acad. Sci. Hungar.,
{\bf{102}} {\rm (1--2)},
123--145,
1989.

[14]
Wikipedia,
{\it{Resultant.}}



\enddocument